\documentclass[a4paper,fleqn]{cas-sc}

\usepackage[authoryear,longnamesfirst]{natbib}

\def\tsc#1{\csdef{#1}{\textsc{\lowercase{#1}}\xspace}}
\tsc{WGM}
\tsc{QE}

\newtheorem{theorem}{Theorem}
\newtheorem{lemma}[theorem]{Lemma}
\newdefinition{rmk}{Remark}
\newproof{Proof}{Proof}

\begin{document}

\let\WriteBookmarks\relax
\def\floatpagepagefraction{1}
\def\textpagefraction{.001}

\shorttitle{}    

\shortauthors{}  

\title [mode = title]{Mass-Conserving Saddle Dynamics via Generalized Inner Product: Theory, Algorithms, and Applications} 


%

\author[1]{Longjing Li}

\affiliation[1]{
            organization={School of Mathematical Sciences, Laboratory of Mathematics and Complex Systems, Ministry of Education,},
            addressline={Beijing Normal University}, 
            city={Beijing},
            postcode={100875}, 
            country={China}
            }
\ead{202321130105@mail.bnu.edu.cn}

\author[1]{Guanghua Ji}

\ead{ghji@bnu.edu.cn}

\cortext[1]{Corresponding author}

\author[2]{Zhen Xu}
\cormark[1]
\affiliation[2]{organization={School of Mathematics, Statistics and Mechanics},
            addressline={Beijing University of Technology}, 
            city={Beijing},
            postcode={100124}, 
            country={China}}
\ead{xuzhenmath@bjut.edu.cn}
\renewcommand{\thefootnote}{\fnsymbol{footnote}}

\begin{abstract}
To reveal the effect of the inner product choice, we present a unified formulation of saddle dynamics for the functional $F$ with a mass constraint under different inner products. We establish the equivalence between the index-$k$ saddle points and the linearly stable steady states of the corresponding dynamics.
Further, we present the dynamics with discrete $H^{-1}$ and $L^2$ inner products and numerically verify the convergence orders of both dynamics.
Finally, we apply the method to a phase field model with driving force under Neumann and periodic boundary conditions. The results uncover previously unreported saddle points and their connectivity, highlighting how the choice of inner product enriches the solution landscape in conservative systems.
\end{abstract}




\maketitle

\begin{abstract}
To reveal the effect of the inner product choice, we present a unified formulation of saddle dynamics for the functional $F$ with a mass constraint under different inner products. We establish the equivalence between the index-$k$ saddle points and the linearly stable steady states of the corresponding dynamics.
Further, we present the dynamics with discrete $H^{-1}$ and $L^2$ inner products and numerically verify the convergence orders of both dynamics.
Finally, we apply the method to a phase field model with driving force under Neumann and periodic boundary conditions. The results uncover previously unreported saddle points and their connectivity, highlighting how the choice of inner product enriches the solution landscape in conservative systems.
\end{abstract}

\section{Introduction}
The phase field method is a computational approach that simulates microstructure evolution using continuous field variables, providing a powerful continuum framework for studying materials science, biological science, and engineering applications \cite{Cahn1958258,phase_field_simulate,phase_field_microstructure}.
Due to non-convexity of the energy functional, the phase field model has multiple critical points. A fundamental problem is to locate all critical points, including both local minima and saddle points. 
From Morse theory \cite{1963Morse}, the solution landscape \cite{Yin2021CGiSD} is introduced to classify nondegenerate critical points with the Morse index, which is the maximal dimension of a subspace on which the Hessian is negative definite. \par

To locate these critical points, various numerical methods have been developed for saddle points, including the iterative minimization formulation \cite{imf,convex_IMF}, minimax-type methods \cite{minmax_Li_zhou,Liu20232361}, and the string method \cite{2002String,simplify_String,String_phase}. For phase field models with mass conservation, numerical methods must incorporate constraints, leading to developments. The constrained string method \cite{constrained_string} and the constrained minimax method \cite{constrain_minimax} introduce the Lagrange multiplier for constraint. The projected iterative minimization \cite{Gu2021} selects the $H^{-1}$ metric to obtain a conservative system and gives a new algorithm in the $L^2$ metric with projection, which is faster than the $H^{-1}$ metric.
However, a systematic comparison of the saddle dynamics and the resulting solution landscape under the $H^{-1}$ and projected $L^2$ inner products, whether they are equivalent or how they differ, remains unexplored.

Within the surface walking methods, there exists dimer-type method \cite{dimer_saddle,dimer_saddle_pre}, gentlest ascent dynamics (GAD) \cite{GAD,GAD_efficient,GAD_nogradient}, and high-index saddle dynamics (HiSD). In particular, Yin et. al. propose high-index saddle dynamics (HiSD) and establish the equivalence between the saddle points of the functional and stable steady states of the saddle dynamics \cite{Yin2019A3576}. 
Subsequent works have extended HiSD to reduce dependence of the initial value \cite{improved_HiSD}, handle systems without gradients \cite{Yin2021CGiSD}, and incorporate linear constraints via constrained HiSD, which defines dynamics on a manifold \cite{2020Constrained}. 
For mass-conserving systems, projection saddle dynamics (PSD), a special case of constrained HiSD, has been applied to the diblock copolymer-homopolymer model with soft confinement \cite{2021Solution}. 
While these above work has primarily analyzed and applied these methods under the inner product and space of $L^2$, they can, in principle, be formulated with other inner products such as $H^{-1}$. 
However, a systematic framework that unifies the saddle dynamics with different inner products remains lacking. Although saddle dynamics with $H^{-1}$ inner product is often considered computationally more expensive, it may reveal distinct solution landscapes to offer more connections and critical points, thereby revealing geometric connectivity patterns that differ from those induced by the standard $L^2$ metric, even though the underlying Morse homology remains invariant \cite{Morse_Homology}.\par 

In this work, we propose a conservative generalized inner product saddle dynamics (CGiSD) that provides a unified formulation for mass-conserving saddle point searching under different inner products. We prove the linear stability of CGiSD under certain assumptions. Additionally, We implement the CGiSD with two specific discrete inner products, projected $L^2$ and $H^{-1}$, for the Ginzburg–Landau functional with a driving force, which serves as a diffuse-interface model for predicting critical nuclei morphology in solid-state phase transformations. We compared the solution landscapes obtained by the two dynamics under periodic and Neumann boundary conditions. 
The main contributions of this article include:
\begin{enumerate}
    \item Theoretically, we establish a unified saddle dynamics framework for mass-conserving systems under different inner products, accompanied by a rigorous linear stability analysis. To the best of our knowledge, this is the first such unified treatment.
    \item Numerically, we implement the dynamics with discrete $L^2$ and $H^{-1}$ inner products for a phase-field model with driving force, and systematically compare the resulting solution landscapes under periodic and Neumann boundary conditions.
    \item Our results reveal that while the two inner-product choices yield identical solution landscapes in certain parameter regimes, they diverge when the interface becomes sharper or the driving force increases, especially under Neumann boundary conditions. These findings highlight the crucial role of inner-product selection in saddle dynamics and provide practical guidance for simulating constrained phase-field systems.   
\end{enumerate}

The remainder of this paper is organized as follows. In Section 2, we derive the CGiSD framework, analyze its theoretical properties, and prove the equivalence of saddle points under different inner products. Section 3 presents numerical experiments, including convergence tests and comparative studies of solution landscapes under various parameters and boundary conditions. Finally, Section 4 concludes the paper with a summary and outlook for future work.


\section{Conserved generalized inner product saddle dynamics}

We now explain the motivation for introducing the conserved generalized inner product saddle dynamics. To minimize a functional \( \mathcal{F}(\phi) \) of function $\phi$ defined  over a domain $\Omega=[0,L]\times[0,L]\subset\mathbb{R}^2$, subject to the constraint of mass conservation 
\begin{equation}
    \int_\Omega \phi(x)\ dx=constant.
    \label{functional_constraint}
\end{equation}
Two classical approaches are the conservative Allen–Cahn equation (conservative AC) \cite{1992Nonlocal,KIM201411} and the Cahn–Hilliard equation (CH) \cite{2000Models}. The conservative $L^2$ gradient flow is
\begin{align}
	\phi_t  =-M_1\operatorname{P}( \frac{\delta \mathcal{F}(\phi)}{\delta \phi} ),\  (x,t) \in \Omega\times (0,T), 
    \label{PAC}
\end{align}
and the $H^{-1}$ gradient flow is
\begin{align}
	\phi_t  =M_2\Delta( \frac{\delta\mathcal{F}(\phi)}{\delta \phi} ),\  (x,t) \in \Omega\times (0,T),
    \label{CH}
\end{align}
where $\frac{\delta \mathcal{F}(\phi)}{\delta \phi}$ is the variational derivative of the functional $\mathcal{F}$ and the projection operator of mass-conservation $\mathrm{P}:L^2(\Omega)\rightarrow L^2(\Omega)$ is defined as $\mathrm{P}(\phi) = \phi - \frac{1}{|\Omega|}\int_{\Omega} \phi\ dx$. $M_1$ and $M_2$ are positive mobilities. Mathematically, the conservative $L^2$ gradient flow can be formulated as a $L^2$-gradient flow with a Lagrange multiplier, while the Cahn-Hilliard equation is a gradient flow in $H^{-1}$, which inherently ensures mass conservation.\par
For numerical calculations, we need to discretize the functional $\mathcal{F}(\phi)$ in space with uniform mesh size $h$ to obtain the discrete functional $F(\boldsymbol{u})$, where $\boldsymbol{u}\in\mathbb{R}^d$ is the vector of grid point values approximating $\phi$. We denote the discrete variational derivative of discrete functional $F$ by $\nabla F(\boldsymbol{u})\in\mathbb{R}^d$.\par
The spatially discrete conservative $L^2$ gradient flow is
\begin{align}
	\boldsymbol{u}_t = -M_1\operatorname{P}_h( \nabla F(\boldsymbol{u})),\  t\in (0,T), 
    \label{D-PAC}
\end{align}
where $M_1$ is positive mobility and $\mathrm{P}_h:\mathbb{R}^d\to\mathbb{R}^d$ is the discrete mass-conservation projection. The spatially discrete $H^{-1}$ gradient flow is
\begin{align}
	\boldsymbol{u}_t = M_2\Delta_h( \nabla F(\boldsymbol{u})),\  t \in (0,T),
    \label{D-CH}
\end{align}
where $M_2$ is positive mobility and $\Delta_h:\mathbb{R}^d \rightarrow \mathbb{R}^d$ denotes the discrete Laplacian under the periodic or homogeneous Neumann boundary condition.\par
Both dynamics share the same steady-state characterization:
\[
\operatorname{P}_h(\nabla F(\boldsymbol{u})) = 0 \quad \Longleftrightarrow \quad \nabla F(\boldsymbol{u}) = \boldsymbol{\mu} \quad \Longleftrightarrow \quad \Delta_h(\nabla F(\boldsymbol{u})) = 0,
\]
where $\boldsymbol{\mu}\in\mathbb{R}^d$ is a constant vector. Nevertheless, the evolutionary paths of the discrete conservative $L^2$ gradient flow \eqref{D-PAC} and $H^{-1}$ gradient flow \eqref{D-CH} are not necessarily identical. \par

\begin{figure}[h]
    \centering
    \includegraphics[scale=0.4]{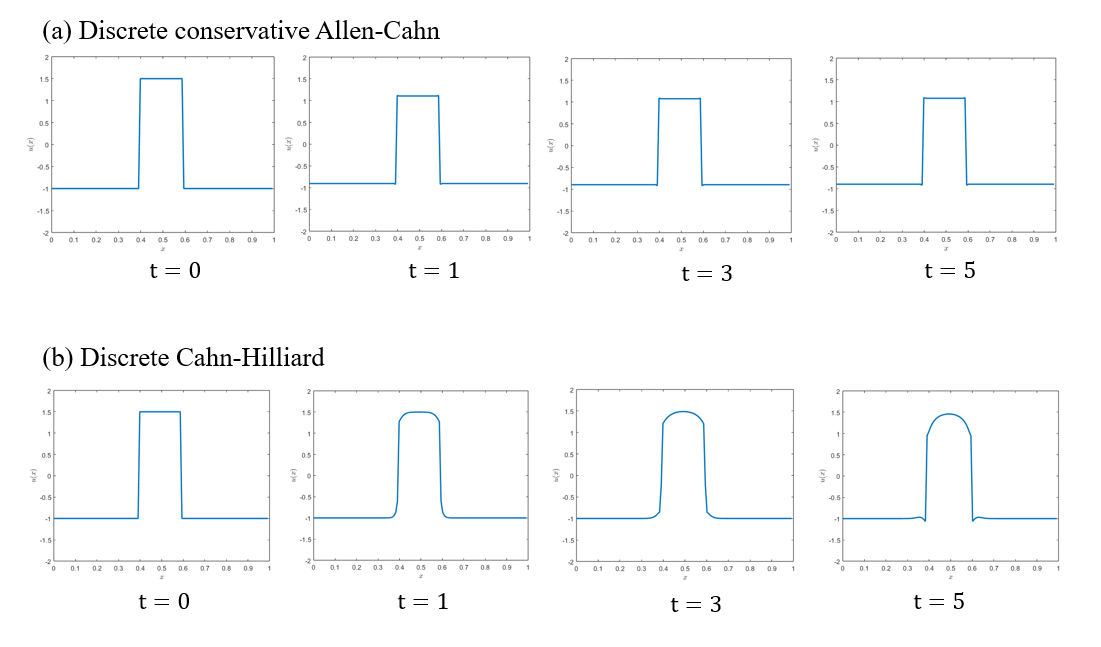}
    \caption{The 1D (a) discrete conservative Allen-Cahn and (b) discrete Cahn-Hilliard of the evolution $\boldsymbol{u}$ with the initial condition $\boldsymbol{u}^0$ and $\delta t=10^{-4}$. Snapshots of the numerical approximation $\boldsymbol{u}$ are taken at $t$ = 0, 1, 5, 10.}
    \label{compare_example1}
\end{figure}

For example, we consider the diffuse-interface functional $\mathcal{F}(\phi)=\int_\Omega (\frac{\epsilon^2}{2}|\nabla \phi|^2+\frac{1}{4}(\phi^2-1)^2)\ dx$, where $\Omega=[0,1]$ and $\epsilon^2=10^{-6}$. We impose the periodic boundary conditions and discretize the domain via the finite difference method on a uniform grid of 128 points. 
Figure~\ref{compare_example1} shows the 1D numerical solutions $\boldsymbol{u}$ of both schemes under the same initial condition $\boldsymbol{u}^0$, whose components are defined as
\[
(\boldsymbol{u}^0)_i=\begin{cases}
	1.5, & 0.4 \leq \frac{i}{128}<0.6 \\
	-1, & else\\
\end{cases}.
\]
The mobilities are chosen as $M_1=1,\, M_2=\frac{1}{128^2}$, which modifies the speed of evolution comparable.
Taking this functional as an example, we then compare the solution landscapes obtained by the two dynamical systems.

It is observed that the discrete AC equation uses the global projection method to preserve mass conservation (Figure \ref{compare_example1} (a)), whereas the discrete CH equation preserves mass conservation achieved through local diffusion (Figure \ref{compare_example1} (b)).
This distinction implies that, for problems with multiple solutions, these two dynamical systems are not necessarily equivalent along the evolution path and may converge to different equilibrium points even from the same initial condition.

\begin{figure}[h]
	\centering
    \includegraphics[scale=0.4]{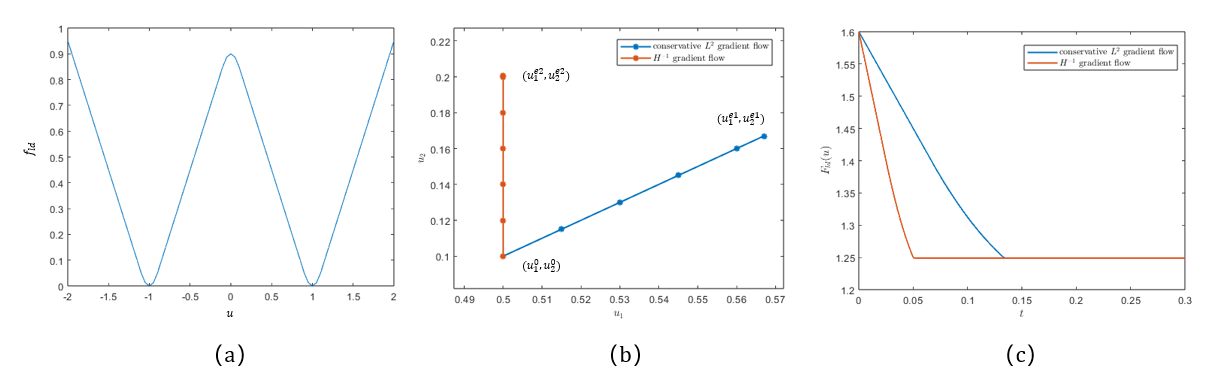}
    \caption{(a) The graph of the potential energy function $F_{dp}\  (\epsilon=0.1)$. (b) The trajectory of $(u_1,u_2)$ from dynamical discrete conservative $L^2$ gradient flow and $H^{-1}$ gradient flow of functional $F_{ld}$. (c) The energy functional $F_{ld}(\boldsymbol{u})$ evolves along the discrete conservative $L^2$ gradient flow and discrete $H^{-1}$ gradient flow.}
    \label{compare_example2}
\end{figure}

To clarify the difference between two dynamical systems, we then consider a modified linear double-well potential function $f_{ld}$, as shown in Figure \ref{compare_example2}(a) and

\begin{equation}
f_{ld}(u)=
\begin{cases}
	-(u+1+\frac{\epsilon_m}{2}), & u<-1-\epsilon_m \\
	\frac{1}{2\epsilon_m}(u+1)^2, & -1-\epsilon_m \leq u < -1+\epsilon_m\\
    u+1-\frac{\epsilon_m}{2}, & -1+\epsilon_m \leq u < -\epsilon_m\\
    -\frac{1}{2\epsilon_m}u^2+1-\epsilon_m, & -\epsilon_m \leq u < \epsilon_m\\
    -u+1-\frac{\epsilon_m}{2}, & \epsilon_m \leq u <1-\epsilon_m\\
    \frac{1}{2\epsilon_m}(u-1)^2, & 1-\epsilon_m \leq u < 1+\epsilon_m\\
    u-1-\frac{\epsilon_m}{2}, & 1+\epsilon_m \leq u
\end{cases}.
\end{equation}
We set the parameters $\epsilon_m=10^{-3}$ and define the functional $F_{ld}(\boldsymbol{u})=\sum_{i=1}^4 f_{ld}(u_i)$, where $\boldsymbol{u}=(u_1,u_2,u_3,u_4)$. Then, the minimum of functional $F_{ld}(\boldsymbol{u})$ with constraint $\sum_{i=1}^4 u_i = 0.6$
can be obtained by evolving the dynamical systems\eqref{D-PAC} and \eqref{D-CH} from initial value $\boldsymbol{u}^0=(0.5,0.1,-0.8,0.8)$. 

The solutions of both discrete conservative $L^2$ gradient flow and $H^{-1}$ gradient flow exist and can be solved analytically. The discrete conservative $L^2$ gradient flow converges to the equilibrium  
\[
\boldsymbol{u}^{\text{e1}} = \left(u_1^{\text{e1}},u_2^{\text{e1}},u_3^{\text{e1}},u_4^{\text{e1}} \right ) = \left(0.5+\frac{2+\epsilon_m}{3},\ 0.1+\frac{2+\epsilon_m}{3},\ -1-\epsilon_m,\ 0.8+\frac{2+\epsilon_m}{3} \right),
\]  
while the discrete $H^{-1}$ one converges to another equilibrium  
\[
\boldsymbol{u}^{e2} =  \left(u_1^{\text{e2}},u_2^{\text{e2}},u_3^{\text{e2}},u_4^{\text{e2}}\right) =\left(0.5,\ 0.2+\frac{\epsilon_m}{2},\ -1-\epsilon_m,\ 0.9+\frac{\epsilon_m}{2} \right).
\] For brevity of visualization, we only show the first two components $(u_1,u_2)$ of the four-dimensional solution vector $\boldsymbol{u}$; the evolution of $(u_3,u_4)$ is analogous and thus omitted.

From Figure \ref{compare_example2} (b), we observe that the solutions $(u_1(t),u_2(t))$ of these two systems evolve with different trajectories and converge to two distinct equilibrium points $(u_1^{e1},u_2^{e1})$ and $(u_1^{e2},u_2^{e2})$, respectively, starting from the same initial value $\boldsymbol{u}^0$. Figure \ref{compare_example2} (c) indicates that the energy $F_{ld}(\boldsymbol{u})$ decreases along the trajectories of both the discrete conservative $L^2$ gradient flow and $H^{-1}$ gradient flow, and becomes stable after a certain period of time. The discrete energy is dissipated along both flows, which is a fundamental feature of gradient flows.\par
Therefore, these two saddle dynamics corresponding to the conservative $L^2$ gradient flow and the $H^{-1}$ gradient flow may follow different evolutionary paths and converge to different saddle points. To explore and analyze the relationship of two saddle dynamics of different inner products, we propose a conservative generalized inner product saddle dynamics (CGiSD).

\subsection{Formulation of conservative generalized inner product saddle dynamics}
Let $\mathcal{H}$ be a Hilbert space of $d$ dimensional and consider a functional $F:\mathcal{H}\rightarrow\mathbb{R}$ that is twice Fr\'{e}chet differentiable. For the Hilbert space $\mathcal{H}=\mathbb{R}^d$ with the standard inner product $\langle \boldsymbol{u},\boldsymbol{v}\rangle=\boldsymbol{u}^\top\boldsymbol{v}$, the derivative and Hessian of functional $F$ are denoted by $\nabla F$ and $\nabla^2F$, respectively. According to Morse theory, a index-$k$ saddle point $\boldsymbol{\hat{u}}$ of functional $F$ in $\mathbb{R}^d$ is $\nabla F(\boldsymbol{\hat{u}}) =0$, and the invertible Hessian  $\nabla^2 F(\boldsymbol{\hat{u}})$ has $k$ negative eigenvalues.\par 
For unconstrained saddle point researching, the high-index saddle dynamics \cite{Yin2019A3576} aims to locate index-$k$ saddle points:\par
\begin{equation}
    \begin{cases}
    	\dot{\boldsymbol{u}}=-\left(I-2\sum_{i=1}^{k}\boldsymbol{v}_{i}\boldsymbol{v}_{i}^{\top}\right)\nabla F(\boldsymbol{u}), \\
    	\dot{\boldsymbol{v}}_{i}=-\left(I-\boldsymbol{v}_{i}\boldsymbol{v}_{i}^{\top}-2\sum_{j=1}^{i-1}\boldsymbol{v}_{j}\boldsymbol{v}_{j}^{\top}\right)\nabla^2F(\boldsymbol{u})[\boldsymbol{v}_{i}], & i=1,\ldots,k.
    \end{cases}
    \label{HIOSD}
\end{equation}

For constrained saddle point researching, we firstly introduce a weighted inner product $\langle\boldsymbol{u},\boldsymbol{v}\rangle_w = \sum_{i=1}^d w_i u_i v_i$,
where the weight $\boldsymbol{w}=(w_1,...,w_d)^\top,\ w_i>0$ and $\sum_{i=1}^d w_i=1$. The continuous functional constraint \eqref{functional_constraint} can be approximated with
\begin{equation}
    \overline{\boldsymbol{u}}(t)=\langle \boldsymbol{u},\boldsymbol{1}\rangle_w=\overline{\boldsymbol{u}}(0),\ \forall\ t>0,
\end{equation}
where $\boldsymbol{1}=(1,...,1)^\top \in \mathbb{R}^d$. For simplicity, the derivative and Hessian of functional $F$ under inner product $\langle \boldsymbol{u},\boldsymbol{v}\rangle_w$ are denoted by $\nabla F$ and $\nabla^2 F$. To satisfy constraint, the variational derivative and the eigenvectors belonging to the space $\mathbb{R}^d_0=\{\boldsymbol{u}\in \mathbb{R}^d: \langle \boldsymbol{1},\boldsymbol{u}\rangle_w = 0\}$ are required. Define an operator $\operatorname{A}:\mathbb{R}^d\rightarrow \mathbb{R}^d$ that is positive semi-definite and self-adjoint whose kernel is $\ker(\operatorname{A})=span\{\boldsymbol{1}\}$. The operator $\operatorname{A}$ has certain properties that will be used in the following proof.

\begin{lemma}
\label{lem:A_property}
Let $\operatorname{A} : \mathbb{R}^d \to \mathbb{R}^d$ be a self-adjoint and positive semi-definite linear operator with kernel $\operatorname{ker}(\operatorname{A}) = \operatorname{span}\{\mathbf{1}\},\text{where}\ \mathbf{1} = (1,1,\dots,1)^T \in \mathbb{R}^d$. Let \(\mathbb{R}^d_0 := \{\mathbf{u} \in \mathbb{R}^d: \langle \boldsymbol{1},\boldsymbol{u} \rangle_w = 0\}\) be the zero-mean subspace. Denote by \(\operatorname{A}_0\) the restriction of $\operatorname{A}$ to \(\mathbb{R}^d_0\). Here $\mathbb{R}^d_0$ and $\mathbb{R}^d$ are equipped with the inner product $\langle \cdot,\cdot\rangle_w$. Then the following properties hold:
\begin{enumerate}
    \item The operator $\operatorname{A}_0 : \mathbb{R}_0^d \to \mathbb{R}_0^d$ is invertible, positive definite, and self-adjoint. Its inverse operator $\operatorname{A}_0^{-1}: \mathbb{R}_0^d \to \mathbb{R}_0^d$ is positive definite and self-adjoint.
    \item Define the projection operator $\mathrm{P}_{\mathrm{h}} : \mathbb{R}^d \to \mathbb{R}^d$, $\mathrm{P}_{\mathrm{h}}(\boldsymbol{u}) = \boldsymbol{u} - \overline{\boldsymbol{u}}\boldsymbol{1}, \ \text{where }   \overline{\boldsymbol{u}}=\langle\boldsymbol{u},\boldsymbol{1}\rangle_w.$
    $\operatorname{P}_h$ is a special case of $\operatorname{A}$ and it is idempotent, i.e., $\mathrm{P}_{\mathrm{h}}^2 = \mathrm{P}_{\mathrm{h}}$. Moreover, for any $\boldsymbol{v} \in \mathbb{R}_0^d$, we have $\mathrm{P}_{\mathrm{h}}(\boldsymbol{v}) = \boldsymbol{v}$.
    \item The following operator identities hold: $\operatorname{A}_0^{-1} \operatorname{A}=\operatorname{P}_h$, $\operatorname{A}\operatorname{A}_0^{-1} =\operatorname{I}$, and $\operatorname{A}_0^{-1}\operatorname{A}_0 =\operatorname{A}_0 \operatorname{A}_0^{-1}=\operatorname{I}$.
\end{enumerate}
\end{lemma}

\begin{Proof}
(1) Since $\ker(\operatorname{A})=span\{\boldsymbol{1}\}$, the restriction $\operatorname{A}_0:\mathbb{R}^d_0\rightarrow\mathbb{R}^d_0$ satisfies $\ker(\operatorname{A}_0)=\{\boldsymbol{0}\}$; thus $\operatorname{A}_0$ is invertible. \\
Define $Im(\operatorname{A}):=\{\operatorname{A}\boldsymbol{u}\,|\,\boldsymbol{u}\in \mathbb{R}^d\}$. From the fundamental subspace theorem, 
\[
Im(\operatorname{A}) = \ker (\operatorname{A}^\top)^\perp=\ker(\operatorname{A})^\perp =\mathbb{R}^d_0.
\]
Next, we prove $Im(\operatorname{A}_0):=\{\operatorname{A}_0\boldsymbol{u}\,|\,\boldsymbol{u}\in \mathbb{R}^d_0\}=\mathbb{R}_0^d$. For $Im(\operatorname{A}_0)\subset \mathbb{R}^d_0$,
\[
\forall\, \boldsymbol{u}\in \mathbb{R}^d_0,\, \operatorname{A}_0\boldsymbol{u}=A\boldsymbol{u}\in Im(\operatorname{A})=\mathbb{R}^d_0.
\]
For $Im(\operatorname{A}_0) \supset \mathbb{R}^d_0$,
\[
\forall\, \boldsymbol{y}\in \mathbb{R}^d_0,\,\exists\, \boldsymbol{u} \in \mathbb{R}^d,\, \operatorname{A} \boldsymbol{u} =\boldsymbol{y}.
\]
From the orthogonal decomposition theorem, there exists $\boldsymbol{u}=\boldsymbol{u}_1+\boldsymbol{u}_2$, $\boldsymbol{u}_1 \in \mathbb{R}^d_0$, $\boldsymbol{u}_2 \in \ker(A)$. Therefore,
\[
\operatorname{A}_0\boldsymbol{u}_1 = \operatorname{A}\boldsymbol{u}_1=\operatorname{A}\boldsymbol{u}=\boldsymbol{y}.
\]
Next, we prove that $A_0$ is positive definite by contradiction. 
$$\forall\, \boldsymbol{u} \in \mathbb{R}^d_0,\ \langle\boldsymbol{u},\operatorname{A}_0\boldsymbol{u}\rangle_w=\langle\boldsymbol{u}, \operatorname{A} \boldsymbol{u}\rangle_w\geq 0,$$
since $\operatorname{A}$ is semi-positive definite. 
Now we assume
$$\exists\  \boldsymbol{u_0}\neq 0,\langle\boldsymbol{u_0},\operatorname{A}_0\boldsymbol{u_0}\rangle_w = 0.$$
Therefore, $\operatorname{A}_0 \boldsymbol{u_0} =0$, which contradicts the invertibility of $\operatorname{A}_0$. $\operatorname{A}_0$ is symmetric by definition since $\operatorname{A}$ is symmetric.\par
(2) $\operatorname{P}_h$ is idempotent since
\[
        \operatorname{P}_h^2(\boldsymbol{u}) = \operatorname{P}_{\mathrm{h}}(\boldsymbol{u} - \overline{\boldsymbol{u}}\boldsymbol{1}) = (\boldsymbol{u} - \overline{\boldsymbol{u}}\boldsymbol{1}) - \overline{(\boldsymbol{u} - \overline{\boldsymbol{u}}\boldsymbol{1})}\boldsymbol{1},
\]
From $\boldsymbol{w}^T \boldsymbol{1}=1$, we have
\[
\overline{(\boldsymbol{u} - \overline{\boldsymbol{u}}\boldsymbol{1})}=\boldsymbol{w}^\top(\boldsymbol{u}-\boldsymbol{w}^\top\boldsymbol{u}\boldsymbol{1})=0.
\]
Therefore, $\operatorname{P}_h^2(\boldsymbol{u})=\boldsymbol{u} - \overline{\boldsymbol{u}}\boldsymbol{1}=\operatorname{P}_h(\boldsymbol{u})$.\par

(3) can be verified through direct calculations.\par

\end{Proof}

For any operator $\operatorname{A}$ and weight vector $\boldsymbol{w}$ satisfying Lemma \ref{lem:A_property}, Hilbert space $\mathbb{R}^d_0$ with inner product $\langle \boldsymbol{u},\boldsymbol{v}\rangle_{A}=\langle \boldsymbol{u}, \operatorname{A}^{-1}_0 \boldsymbol{v}\rangle_w$ is well-defined. For the functional $F$, the inner product $\langle\cdot,\cdot\rangle_A$ on the tangent space $\mathbb{R}^d_0$ of the manifold $\boldsymbol{u}+\mathbb{R}^d_0$ introduces the gradient $\nabla_A F(\boldsymbol{u})$
\[
    \langle\nabla_A F(\boldsymbol{u}),\boldsymbol{\psi_1}\rangle_A =\lim\limits_{\epsilon\rightarrow0} \frac{F(\boldsymbol{u}+\epsilon \boldsymbol{\psi}_1)-F(\boldsymbol{u})}{\epsilon}, \forall\, \boldsymbol{u} \in \mathbb{R}^d,\, \boldsymbol{\psi}_1\in \mathbb{R}^d_0
\]
and the Hessian $\nabla^2_A F(\boldsymbol{u})$
\[
    \langle\nabla^2_A F(\boldsymbol{u})[\boldsymbol{\psi}_1],\boldsymbol{\psi}_2\rangle_A =\lim\limits_{\epsilon\rightarrow0} \frac{\langle \nabla_AF(\boldsymbol{u}+\epsilon \boldsymbol{\psi}_2)-\nabla_AF(\boldsymbol{u}),\boldsymbol{\psi}_1\rangle_A}{\epsilon},\forall\, \boldsymbol{u} \in \mathbb{R}^d,\boldsymbol{\psi}_1,\boldsymbol{\psi}_2\in \mathbb{R}^d_0.\]
By the definition of the induced gradient and Hessian, we have the relations $\nabla_A F(\boldsymbol{u})=\operatorname{A}\nabla F(\boldsymbol{u})$ and $\nabla^2_AF(\boldsymbol{u})=\operatorname{A} \nabla^2 F(\boldsymbol{u})$.\par

The construction of conservative generalized inner product saddle dynamics consists of two steps. The first step is to construct the evolution of $\boldsymbol{u}$ under the inner product $\langle\cdot,\cdot\rangle_A$, and the second step is to compute the eigenvectors corresponding to the $k$ smallest eigenvalues.\par

Let $\boldsymbol{\hat{u}}$ be an index-$k$ saddle point, i.e., the gradient $\nabla_A F(\hat{\boldsymbol{u}})=0$ and Hessian $\nabla^2_AF(\boldsymbol{\hat{u}})$ is invertible with $k$ negative eigenvalues $\hat{\lambda}_1,...,\hat{\lambda}_k$. From the following Lemma \ref{A_orthogonal}, there exists a set of eigenvectors $\hat{\boldsymbol{v}}_1,...,\hat{\boldsymbol{v}}_k\in\mathbb{R}^d_0$ of $\nabla^2_A F(\hat{\boldsymbol{u}})$ corresponding to the eigenvalues $\hat{\lambda}_1,...,\hat{\lambda}_k$, satisfying $\langle \hat{\boldsymbol{v}}_i,\hat{\boldsymbol{v}}_j\rangle_A=\delta_{ij}$.
\begin{lemma} 
 If $F:\mathbb{R}^d\rightarrow\mathbb{R}$ is a $\mathcal{C}^3$ functional, then for every $\boldsymbol{u} \in \mathbb{R}^d$, the operator $\nabla^2_AF(\boldsymbol{u})$ admits a set of eigenvectors $\{\boldsymbol{v}_1,...,\boldsymbol{v}_{d-1}\} \subset \mathbb{R}^d_0$ that are orthogonal with respect to the inner product $\operatorname{A}$, i.e., $\langle \boldsymbol{v}_i,\boldsymbol{v}_j\rangle_A=\delta_{ij}$.
 \label{A_orthogonal}
\end{lemma}
\begin{Proof}
$\ \forall\ \boldsymbol{u} \in \mathbb{R}^d$ , if $ \boldsymbol{v} \in \mathbb{R}^d_0$ is a eigenvector of 
$\nabla^2_AF(\boldsymbol{u})=\operatorname{A}\nabla^2F(\boldsymbol{u})$,
\[ \operatorname{A} \nabla^2 F(\boldsymbol{u})\boldsymbol{[v]}=\lambda \boldsymbol{v}\Leftrightarrow \operatorname{P}_h \nabla^2F(\boldsymbol{u})\operatorname{P}_h[\boldsymbol{v}]=\lambda \operatorname{A}^{-1}_0 \boldsymbol{v}\  ,\ \lambda\in \mathbb{R}
\]
$\operatorname{P}_h \nabla^2F(u)\operatorname{P}_h[\boldsymbol{v}]=\lambda \operatorname{A}^{-1}_0 \boldsymbol{v}$ is generalized eigenvalue problem, the conclusion hold from $\operatorname{P}_h\nabla^2F(\boldsymbol{u})\operatorname{P}_h$ is self-adjoint and $\operatorname{A}^{-1}_0$ is positive definite.
\end{Proof} \par

Similar to high-index saddle dynamics \cite{Yin2019A3576}, by setting $\hat{\mathcal{V}}=span\{\hat{\boldsymbol{v}}_1,...,\hat{\boldsymbol{v}}_k\}$, then $\boldsymbol{\hat{u}}$ is a local maximum in the $k$-dimensional manifold $\hat{\boldsymbol{u}}+\hat{\mathcal{V}}$ and a local minimum in $\hat{\boldsymbol{u}}+\hat{\mathcal{V}}^\perp$. We assume that $\boldsymbol{u}$ is sufficiently close to the index-$k$ saddle point $\hat{\boldsymbol{u}}$ and the Hessian $\nabla^2_AF(\boldsymbol{u})$ also has $k$ negative eigenvalues. To make $\hat{\boldsymbol{u}}$ approach with $\boldsymbol{u}$, $\boldsymbol{u}$ should ascend in manifold $\hat{\boldsymbol{u}}+\hat{\mathcal{V}}$ and descend in $\hat{\boldsymbol{u}}+\hat{\mathcal{V}}^\perp$. 
Approximating $\hat{\boldsymbol{u}}+\hat{\mathcal{V}}$ by $\boldsymbol{u}+\mathcal{V}$, we obtain the dynamical system of $\boldsymbol{u}$ with inner product $\langle\cdot,\cdot\rangle_A$ as follows:
\begin{align*}
	\dot{\boldsymbol{u}} &= 
    \sum_{i=1}^k\langle\nabla_AF(\boldsymbol{u}),\boldsymbol{v}_i\rangle_A \boldsymbol{v}_i+(-\nabla_AF(\boldsymbol{u})+\sum_{i=1}^k\langle\nabla_AF(\boldsymbol{u}),\boldsymbol{v}_i\rangle_A\boldsymbol{v}_i)\\
    &=-\nabla_AF(\boldsymbol{u}) +2\sum_{i=1}^k\langle \nabla_AF(\boldsymbol{u}),\boldsymbol{v}_i\rangle_A \boldsymbol{v}_i\\
    &= -(I-2\sum_{i=1}^k\boldsymbol{v}_i(\boldsymbol{w}\circ\boldsymbol{v}_i)^T\operatorname{A}^{-1}_0)\nabla_AF(\boldsymbol{u}),
\end{align*}
where Hadamard product (or entrywise product) of two vectors $\boldsymbol{w},\boldsymbol{v}\in \mathbb{R}^d$ is $\boldsymbol{w}\circ\boldsymbol{v_i}=(w_1v_1,w_2v_2,...,w_dv_d)^T$.

To derive the eigenvectors $\boldsymbol{v}_i\ (1\leq i\leq k)$ of $\operatorname{A}\nabla^2F(\boldsymbol{u})$ , we can transform $\operatorname{A}\nabla^2F(\boldsymbol{u}) [\boldsymbol{v}_i]=\lambda \boldsymbol{v}_i$ into the generalized eigenvector problem $$\operatorname{P}_h\nabla^2F(\boldsymbol{u}) [\operatorname{P}_h \boldsymbol{v}_i]=\lambda_i \operatorname{A}^{-1}_0\boldsymbol{v_i},\quad 1\leq i\leq k,$$
where $\lambda_i$ is the $k$-th smallest eigenvalue of $\nabla^2_AF$. The subspace $\mathcal{V}$ is constructed simultaneously by solving $k$ the generalized Rayleigh quotient with constraints,
\begin{equation}
    \min_{\boldsymbol{v}_i} \frac{\langle \boldsymbol{v}_i, \operatorname{P}_h\nabla^2F(\boldsymbol{u})[\operatorname{P}_h \boldsymbol{v}_i] \rangle_w}{\langle \boldsymbol{v}_i,\boldsymbol{v}_i\rangle_A}  \quad \text{s.t.} \quad \langle \boldsymbol{v}_i, \boldsymbol{v}_j \rangle_A = \delta_{ij}.
\label{raylei-quotient}
\end{equation}
Then we minimize the $k$ Rayleigh quotients (\ref{raylei-quotient}) simultaneously by using the dynamics of $\boldsymbol{v}_i$ as follows
\begin{equation}
    \dot{\boldsymbol{v}}_{i}=-\left(I-\boldsymbol{v}_{i}\boldsymbol{v}_i^{w,\top}\operatorname{A}^{-1}_0-2\sum_{j=1}^{i-1}\boldsymbol{v}_{j}\boldsymbol{v}^{w,\top}_j\operatorname{A}^{-1}_0\right)\nabla^2_AF(\boldsymbol{u})[\operatorname{P}_h\boldsymbol{v}_{i}], i=1,\ldots,k,
    \label{v_d}
\end{equation}
where $\boldsymbol{v}^{w,\top}_i = (\boldsymbol{w} \circ\boldsymbol{v}_i)^\top$. \par

To derive (\ref{v_d}), the Lagrangian function $\mathcal{L}_i (\boldsymbol{v}_i; \xi_1, \ldots, \xi_i)$ of (\ref{raylei-quotient}) is
\begin{equation}
     \mathcal{L}_i= \langle \boldsymbol{v}_i, \operatorname{P}_h\nabla^2F(\boldsymbol{u})[\operatorname{P}_h \boldsymbol{v}_i] \rangle_w - \xi_i (\langle \boldsymbol{v}_i, \boldsymbol{v}_i \rangle_A - 1) - \sum_{j=1}^{i-1} \xi_j \langle \boldsymbol{v}_j, \boldsymbol{v}_i \rangle_A,
    \label{eq:lagrangian}
\end{equation}
where $\xi_1, \ldots, \xi_i$ are Lagrangian multipliers, and the gradient of \eqref{eq:lagrangian} with the inner product $\langle\cdot,\cdot\rangle_A$ is
\begin{equation}
    \frac{\partial}{\partial \boldsymbol{v}_i} \mathcal{L}_i (\boldsymbol{v}_i; \xi_1, \ldots, \xi_{i-1}, \xi_i) = 2\nabla^2_AF(\boldsymbol{u})\boldsymbol{v}_i - 2\xi_i \boldsymbol{v}_i - \sum_{j=1}^{i-1} \xi_j \boldsymbol{v}_j.
    \label{eq:gradient}
\end{equation}
The dynamics of $\boldsymbol{v}_i$ is
\begin{equation}
    \dot{\boldsymbol{v}}_i = -2\nabla^2_AF(u)[\boldsymbol{v}_i] + 2\xi_i \boldsymbol{v}_i + \sum_{j=1}^{i-1} \xi_j \boldsymbol{v}_j.
    \label{v_d_eps}
\end{equation}
Since $\langle \dot{\boldsymbol{v}}_i,\boldsymbol{v}_i\rangle_A =0$ and $\langle \boldsymbol{v}_i,\boldsymbol{v}_i\rangle_A =1$, we take the inner product of both sides of \eqref{v_d_eps} under $\langle \cdot, \cdot \rangle_A $ and obtain
\[
    \xi_i=\langle \nabla^2_AF(\boldsymbol{u})[\operatorname{P}_h\boldsymbol{v}_i],\boldsymbol{v}_i\rangle_A.
\]
To determine $\xi_j$ for $j<i$, we use the differentiated orthogonality condition $\langle \dot{\boldsymbol{v}}_i,\boldsymbol{v}_j\rangle_A+\langle \boldsymbol{v}_i,\dot{\boldsymbol{v}}_j\rangle_A=0$.
Substituting the expressions for $\boldsymbol{v}_i$ and $\boldsymbol{v}_j$ into this relation and using the orthogonality constraints yields
\[
 \quad \xi_j=2\langle\nabla^2_AF(\boldsymbol{u})[\operatorname{P}_h\boldsymbol{v}_i],\boldsymbol{v}_j\rangle_A.
\]

We thus obtain the following conservative generalized inner product saddle dynamics (CGiSD):
\begin{equation}
    \begin{cases}
    	\dot{\boldsymbol{u}}=-\left(I-2\sum_{i=1}^{k}\boldsymbol{v}_{i}\boldsymbol{v}_i^{w,\top}\operatorname{A}^{-1}_0\right)\nabla_A F(\boldsymbol{u}), \\
    	\dot{\boldsymbol{v}}_{i}=-\left(I-\boldsymbol{v}_{i}\boldsymbol{v}_i^{w,\top}\operatorname{A}^{-1}_0-2\sum_{j=1}^{i-1}\boldsymbol{v}_{j}\boldsymbol{v}^{w,\top}_j\operatorname{A}^{-1}_0\right)\nabla^2_AF(\boldsymbol{u})[\operatorname{P}_h\boldsymbol{v}_{i}],\ 1\leq i \leq k.
    \end{cases}
    \label{GiSD}
\end{equation}
\\
\textbf{Remark 1}\\
Let $\operatorname{A}=\operatorname{P}_h$, then $\operatorname{A}^{-1}_0=\operatorname{I}$. Therefore, the projection saddle dynamics method (PSD) \cite{2021Solution} is a special case of CGiSD.

\subsection{Linear stability analysis of CGiSD}
To achieve better stability, a modified term with $\mu >0$ is attached to the dynamics of $\boldsymbol{u}$ and extend $\nabla^2_AF(\boldsymbol{u})[\operatorname{P}_h\boldsymbol{v}_i]$ to reinforce the mass conservation, leading to a modified CGiSD for index-$k$ saddle points,
\begin{equation}
    \begin{cases}
    	\dot{\boldsymbol{u}}=-\left(I-2\sum_{i=1}^{k}\boldsymbol{v}_{i}\boldsymbol{v}_i^{w,\top}\operatorname{A}^{-1}_0\right)\nabla_A F(\boldsymbol{u})-\mu (\overline{\boldsymbol{u}}-\overline{\boldsymbol{u}}(0))\boldsymbol{1}^\top, \\
    	\dot{\boldsymbol{v}}_{i}=-\left(I-\boldsymbol{v}_{i}\boldsymbol{v}_i^{w,\top}\operatorname{A}^{-1}_0-2\sum_{j=1}^{i-1}\boldsymbol{v}_{j}\boldsymbol{v}^{w,\top}_j\operatorname{A}^{-1}_0\right)\nabla^2_AF(\boldsymbol{u})[\operatorname{P}_h\boldsymbol{v}_{i}], & i=1,\ldots,k.
    \end{cases}
    \label{M-GiSD}
\end{equation}
The term $\mu(\overline{\boldsymbol{u}}-\overline{\boldsymbol{u}}(0))\boldsymbol{1}$ drives $\overline{\boldsymbol{u}}$ toward its initial value $\overline{\boldsymbol{u}}(0)$. We note that the original CGiSD already satisfies the mass conservation constraint $\overline{\boldsymbol{u}}(t)=\overline{\boldsymbol{u}}(0)$ for all $t$.
Consequently, the additional term vanishes identically along any trajectory, so the modified CGiSD is equivalent to the original CGiSD.\par
Next, we will prove that the stable states of \eqref{M-GiSD} are exactly $k$-saddle points. To prove this, we need the following Lemma \ref{FA_decomposed}.

\begin{sloppypar}
\begin{lemma} 
Let \(\operatorname{A}\) satisfy the conditions in Lemma~\ref{lem:A_property}, and let \(\nabla^2_AF(\boldsymbol{u}) = \operatorname{A} \nabla^2 F(\boldsymbol{u})\). Suppose \(\boldsymbol{v}_1,\dots,\boldsymbol{v}_{d-1} \in \mathbb{R}^d_0\) are eigenvectors of the operator \(\nabla^2_AF(\boldsymbol{u}) \operatorname{P}_h\) with corresponding eigenvalues \(\lambda_1,\dots,\lambda_{d-1}\). Then the following spectral decomposition holds:
\[
\nabla^2_AF(\boldsymbol{u}) \operatorname{P}_h = \sum_{i=1}^{d-1} \lambda_i \, \boldsymbol{v}_i \, \boldsymbol{v}_i^{w,\top} \operatorname{A}_0^{-1}. 
\] \label{FA_decomposed}
\end{lemma}
\end{sloppypar}
\begin{Proof}
From the Lemma \ref{lem:A_property} (3), we have $\nabla^2_AF(\boldsymbol{u}) = \operatorname{A}_0 \operatorname{P}_h \nabla^2F(\boldsymbol{u})$, then 

\[ \nabla^2_AF(\boldsymbol{u})[\operatorname{P}_h\boldsymbol{v}_i]=\lambda_i \boldsymbol{v}_i\Leftrightarrow \operatorname{A}_0 \operatorname{P}_h \nabla^2F(\boldsymbol{u})[\operatorname{P}_h \boldsymbol{v}_i]=\lambda_i \boldsymbol{v}_i ,\ 1\leq i \leq d-1.
\]
Since $\operatorname{A}_0$ is positive definite, we decompose it as $\operatorname{A}_0=(\operatorname{A}_0^{\frac{1}{2}})^2$ with $\operatorname{A}_0^\frac{1}{2}$ positive definite and symmetric. 
Introducing $\boldsymbol{t}_i =(\operatorname{A}_0^{\frac{1}{2}})^{-1}\boldsymbol{v}_i$ and multiplying the equation by $(\operatorname{A}_0^{\frac{1}{2}})^{-1}\boldsymbol{v}_i$ on the left give
\[\operatorname{A}_0^{\frac{1}{2}} \operatorname{P}_h \nabla^2F(\boldsymbol{u})\operatorname{P}_h\operatorname{A}_0^{\frac{1}{2}}\boldsymbol{t}_i=\lambda_i \boldsymbol{t}_i.
\]

The operator on the left is symmetric under the weighted inner product, hence by the spectral theorem it can be diagonalized:
\[
\operatorname{A}_0^{\frac{1}{2}} \operatorname{P}_h \nabla^2F(\boldsymbol{u})\operatorname{P}_h\operatorname{A}_0^{\frac{1}{2}}=\sum_{i=1}^{d-1}\lambda_i \boldsymbol{t}_i(\boldsymbol{w}\circ\boldsymbol{t}_i)^\top =(\operatorname{A}_0^{\frac{1}{2}})^{-1}\sum_{i=1}^{d-1}\lambda_i \boldsymbol{v}_i(\boldsymbol{w}\circ\boldsymbol{v}_i)^\top(\operatorname{A}_0^{\frac{1}{2}})^{-1}.
\]

Multiplying on the left by $\operatorname{A}_0^{\frac{1}{2}}$, then right-multiplying by $\operatorname{A}_0^{-1}$ yields
\[
\operatorname{A}_0 \operatorname{P}_h \nabla^2F(\boldsymbol{u})\operatorname{P}_h=\sum_{i=1}^{d-1}\lambda_i \boldsymbol{v}_i(\boldsymbol{w}\circ\boldsymbol{v}_i)^\top\operatorname{A}_0^{-1}.
\]
Since $\nabla^2_AF(\boldsymbol{u})\operatorname{P}_h = \operatorname{A}_0 \operatorname{P}_h \nabla^2F(\boldsymbol{u}) \operatorname{P}_h$, the conclusion is proved.
\end{Proof}

\begin{theorem}
    \label{theorem:linear}
	Assume that $F(\boldsymbol{u})$ is a $\mathcal{C}^3$ functional, $\boldsymbol{u}^*\in \mathbb{R}^d$ and the unit vectors $\{\boldsymbol{v}^*_i\}_{i=1}^k\in \mathbb{R}^d_0$. Suppose $\operatorname{Hess}F(\boldsymbol{u}^*)$ is nondegenerate, with eigenvalues \(\lambda_1^* \leq \ldots \leq \lambda_k^*< 0 < \lambda_{k+1}^* \leq \ldots \leq \lambda_{d-1}^*\). 
	Then the following statements are equivalent:
	\begin{enumerate}
		\item[(A)] \((\boldsymbol{u}^*, \boldsymbol{v}_1^*, \ldots, \boldsymbol{v}_k^*)\) is a linearly stable steady state of CGiSD(\ref{M-GiSD}).
		\item[(B)] \(\boldsymbol{u}^*\) is an index-$k$ saddle point of F, \(\{\boldsymbol{v}_i^*\}_{i=1}^k\) are the eigenvectors of $\nabla^2_AF(\boldsymbol{u}^*)$, corresponding to the eigenvalues \(\{\lambda_i^*\}_{i=1}^k\).
	\end{enumerate}
\end{theorem}
\begin{Proof}
	To determine the stable state of CGiSD(\ref{GiSD}), we calculate its Jacobian operator as follows.
	\begin{equation}
		J = \frac{\partial (\dot{\boldsymbol{u}},  \dot{\boldsymbol{v}}_1, \dot{\boldsymbol{v}}_2, \ldots, \dot{\boldsymbol{v}}_k,)}{\partial (\boldsymbol{u},  \boldsymbol{v}_1, \boldsymbol{v}_2, \ldots \boldsymbol{v}_k)} = 
		\begin{bmatrix}
			 J_{\boldsymbol{u}} & J_{\boldsymbol{u}1} & J_{\boldsymbol{u}2} & \cdots & J_{\boldsymbol{u}k}\\
			 * & J_{1} & 0 & \cdots & 0\\
			 * & J_{21} & J_{2} & \cdots & 0 \\
			 \vdots & \vdots & \vdots & \ddots & \vdots\\
			 * & J_{k1} & J_{k2} & \cdots & J_{k}  \\
		\end{bmatrix}.
	\end{equation}
    Certain blocks are presented in detail as follows:
        \begin{align*}
		J_{\boldsymbol{u}} &= \frac{\partial \dot{\boldsymbol{u}}}{\partial \boldsymbol{u}}= -(\operatorname{I}-2\sum_{i=1}^{k}\boldsymbol{v}_i\boldsymbol{v}_i^{w,\top}\operatorname{A}_0^{-1})\nabla^2_AF(\boldsymbol{u})-\mu \boldsymbol{1} \boldsymbol{1}^{w,\top},\\
		J_{\boldsymbol{u}i} &= \frac{\partial \dot{\boldsymbol{u}}}{\partial \boldsymbol{v}_i} = -2\left( \boldsymbol{v}_i^{w,\top} \operatorname{A}^{-1}_0 \nabla_A F(\boldsymbol{u}) I + \boldsymbol{v}_i (\boldsymbol{w}\circ\nabla _AF(\boldsymbol{u}))^\top \operatorname{A^{-\top}_0} \right)\ ,\\
		J_{i} &=\frac{\partial \dot{\boldsymbol{v}_i}}{\partial \boldsymbol{v}_i} = -(I-2\sum_{j=1}^{i}\boldsymbol{v}_j\boldsymbol{v}_j^{w,\top} \operatorname{A}^{-1}_0)\nabla^2_AF(\boldsymbol{u})\operatorname{P}_h+I\boldsymbol{v}_i^{w,\top} \operatorname{A}^{-1}_0\nabla^2_AF(\boldsymbol{u})[\operatorname{P}_h\boldsymbol{v}_i],
	\end{align*}
    where $1\leq i \leq k$ and $\boldsymbol{v}_i^{w,\top}=(\boldsymbol{w}\circ\boldsymbol{v}_i)^\top$.
    $*$ entries in the Jacobian matrix denotes blocks that can be inferred from symmetry or are not needed for the stability analysis.\\
	\textbf{(A)}$\Leftarrow$\textbf{(B)}: We suppose that $(\boldsymbol{u}^*,\boldsymbol{v}_1^*,...,\boldsymbol{v}_k^*)$ is the index-$k$ saddle point of $F$ with mass conservation, which implies $\nabla_AF(\boldsymbol{u}^*)=0$, $\nabla^2_AF(\boldsymbol{u}^*)[\boldsymbol{v}_i^*]=\lambda_i^* \boldsymbol{v}_i^*(1\leq i \leq k)$, $\lambda_i^*<0,\ \boldsymbol{v}_i^*\in\mathbb{R}^d_0$.\par
	Due to $\nabla_AF(\boldsymbol{u}^*)=0$, $J_{\boldsymbol{u}i}^*=0$. So, $J^*$ is partitioned lower triangular matrix.
	The spectrum of $J^*$ is completely determined by matrix $J_{\boldsymbol{u}}^*,J_{1}^*,...,J_{k}^*$. Next, we will prove that $J_{\boldsymbol{u}}^*$ and $J_{i}^*(1\leq i\leq k)$ both have $d$ negative eigenvalues.
	\begin{equation}
		J_{\boldsymbol{u}}^*[\boldsymbol{v}_j^*] = 
		\begin{cases} 
			\lambda_j^* \boldsymbol{v}_j^*, & 1 \leq j \leq i, \\
			-\lambda_j^* \boldsymbol{v}_j^*, & i < j \leq d-1,
		\end{cases}
	\end{equation}
	indicates $J_{\boldsymbol{u}}^*$ has $d-1$ negative eigenvalues $\lambda_1^*,...,\lambda_k^*,-\lambda_{k+1}^*,...,-\lambda_{d-1}^*$. And $\boldsymbol{1}$ is also a eigenvector of $J^*_{\boldsymbol{u}}$ corresponding eigenvalue $-\mu$.\\
    
    We note that $\boldsymbol{v}_j^*\ (\forall 1\leq j\leq d-1)$ are eigenvectors of $J_{i}^*$:
	\begin{equation}
		J_{i}^*[\boldsymbol{v}_j^*] = 
		\begin{cases} 
			(\lambda_j^* + \lambda_i^*) \boldsymbol{v}_j^*, & 1 \leq j \leq i, \\
			(-\lambda_j^* + \lambda_i^*)\boldsymbol{v}_j^*, & i < j \leq d-1,
		\end{cases}
	\end{equation}
\noindent
which indicates that $J_{i}^*$ has $d-1$ negative eigenvalues $\lambda_1^*+\lambda_i^*,...,\lambda_i^*+\lambda_i^*,-\lambda_{i+1}^*+\lambda_{i}^*,...,-\lambda_{d-1}^*+\lambda_{i}^*\ (\forall i =1,...,k)$. And $\boldsymbol{1}$ is also an eigenvector of $J^*_{\boldsymbol{u}}$ corresponding to the eigenvalue $-\mu$. 
Therefore, \((\boldsymbol{u}^*, \boldsymbol{v}_1^*, \ldots, \boldsymbol{v}_k^*)\) is linearly stable for the dynamical system (\ref{GiSD}).\par 

	\textbf{(A)}$\Rightarrow$\textbf{(B)}: Suppose $(\boldsymbol{u}^*,\boldsymbol{v}_1^*,...,\boldsymbol{v}_k^*)$ is a linearly stable steady state, i.e., $\dot{\boldsymbol{u}}^*=0$ and $\dot{\boldsymbol{v}}_i^*=0\ (1\leq i \leq k)$.\par 
	From $\dot{\boldsymbol{v}}_i^*=0$, we have
	\begin{displaymath}
		\left( I-\boldsymbol{v}_i^*(\boldsymbol{w}\circ\boldsymbol{v}_i^*)^{\top}\operatorname{A}^{-1}_0-2\sum_{j=1}^{i-1}\boldsymbol{v}_j^*(\boldsymbol{w}\circ\boldsymbol{v}_j^*)^{\top}\operatorname{A}^{-1}_0)\nabla^2_AF(\boldsymbol{u}^*)\right)[\operatorname{P}_h\boldsymbol{v}_i^*]=0.
	\end{displaymath}
	We now prove the following property by mathematical induction that for $i=1,...k$:
	\begin{equation}
		\label{eigenvalue-prove}
		\nabla^2_AF(\boldsymbol{u}^*)[\operatorname{P}_h \boldsymbol{v}_i^*]=\mu_i^* \boldsymbol{v}_i^*\neq 0,\,\langle \boldsymbol{v}_j^*,\boldsymbol{v}_i^*\rangle_A=0\ (\forall j=1,...,i).
	\end{equation}
	For $i=1$, it's obvious that $\nabla^2_AF(\boldsymbol{u}^*)[\operatorname{P}_h\boldsymbol{v}_1^*]=\mu_1 \boldsymbol{v}_1^*$ and $\mu_1=\langle \nabla^2_AF(\boldsymbol{u}^*)[\boldsymbol{v}_1^*], \boldsymbol{v}_1^* \rangle_A\neq 0$. Assuming that (\ref{eigenvalue-prove}) holds for $1 \leq i \leq l-1$, by taking $i=l$, we have
	\begin{displaymath}
		\left(\nabla^2_A F(\boldsymbol{u}^*)\operatorname{P}_h - \sum_{j=1}^{l-1} 2\mu_j^* \boldsymbol{v}_j^{*} (\boldsymbol{w}\circ\boldsymbol{v}_j^*)^{\top}\operatorname{A}^{-1}_0 \right) \boldsymbol{v}_l^* = \mu_l^* \boldsymbol{v}_l^*.
	\end{displaymath}
	Since $\boldsymbol{v}_1^*,...,\boldsymbol{v}_{l-1}^*$ are eigenvectors of $\nabla^2_A F(\boldsymbol{u}^*)$, the eigenvectors of $\nabla^2_AF(\boldsymbol{u}^*)$ and $(\nabla^2_AF(\boldsymbol{u}^*)-\sum_{j=1}^{l-1}2\mu_j^*\boldsymbol{v}_j^*(\boldsymbol{w}\circ\boldsymbol{v}_j^*)^{\top}\operatorname{A}^{-1}_0)$ are same from Lemma \ref{FA_decomposed}. We then derive
	\begin{displaymath}
		\nabla^2_A F(\boldsymbol{u}^*) [\boldsymbol{v}_l^*] = \mu_l^* \boldsymbol{v}_l^*.
	\end{displaymath}
	Furthermore, we have $\sum_{i=1}^{l-1}\mu_j^* \langle \boldsymbol{v}_j^*,\boldsymbol{v}_l^* \rangle_A=0$. Since the eigenvectors $\{\boldsymbol{v}_i^*\}_{i=1}^l$ are orthogonal to each other and $\mu_l^* \neq 0$, we have $\langle \boldsymbol{v}_i^*,\boldsymbol{v}_j^* \rangle_A =0$. Thus the property (\ref{eigenvalue-prove}) is proved.\par
	From $\dot{\boldsymbol{u}}^*=0$, we have $(I-2\sum_{i=1}^{k}\boldsymbol{v}_i^*(\boldsymbol{w}\circ\boldsymbol{v}_i^*)^{\top}\operatorname{A}^{-1}_0)\nabla_AF(\boldsymbol{u}^*)=0$. From $(I-\sum_{i=1}^{k}\boldsymbol{v}_i^*(\boldsymbol{w}\circ\boldsymbol{v}_i^{*})^{\top}\operatorname{A}^{-1}_0)^2=I$, we have $\nabla_AF(\boldsymbol{u}^*)=0$ which indicates that $\boldsymbol{u}^*$ is a critical point of $F(\boldsymbol{\boldsymbol{u}})$.
	Similarly, $J^*$ can be viewed as a partitioned lower triangular matrix. $J^*_{\boldsymbol{u}}$ and $J^*_i\ (1\leq i \leq k)$ only have negative eigenvalues.\par
	By property (\ref{eigenvalue-prove}), the eigenvectors $\boldsymbol{v}_1^*,...,\boldsymbol{v}_{d-1}^*$ of $J_{\boldsymbol{u}}^*$ correspond to eigenvalues $\{\mu_1^*,...,\mu_k^*,-\mu_{k+1}^*,...,-\mu_{d-1}^*\}$. After repeating the calculations and transformations, $J^*_i\ (1\leq i \leq k)$ only have negative eigenvalues, i.e., $\mu_i<0,i\leq k$ and $\mu_i>0,k+1 \leq i \leq d-1$.\par
Similarly, for each $i$ ($1\le i\le k$), the eigenvectors $\boldsymbol{v}_1^*,\dots,\boldsymbol{v}_{d-1}^*$ of $J_i^*$ correspond to eigenvalues
\[
\mu_1^*+\mu_i^*,\ \dots,\ \mu_i^*+\mu_i^*,\ -\mu_{i+1}^*+\mu_i^*,\ \dots,\ -\mu_{d-1}^*+\mu_i^*,
\] all of which are negative because $\mu_i^*<0$ and $\mu_j^*>0$ for $j>i$. 
Consequently, all diagonal blocks have only negative eigenvalues, implying that $(\boldsymbol{u}^*,\boldsymbol{v}_1^*,\dots,\boldsymbol{v}_k^*)$ is indeed a linearly stable steady state. Moreover, the inequalities $\mu_1^* < \mu_2^* < \dots < \mu_k^* < 0 \le \mu_{k+1}^* \le \dots \le \mu_{d-1}^*$ hold.
\end{Proof}

\subsection{The equilibrium point equivalence of CGiSD}
In this section, we will prove that the steady stable states of two saddle dynamics systems generated by $\operatorname{L}$ and $\operatorname{T}$ are equivalent, where $\operatorname{L}$ and $\operatorname{T}$ are two different semi-positive definite and self-adjoint operator with $\ker(\operatorname{L})=\ker(\operatorname{T})=span\{\boldsymbol{1}\}$.   
\begin{theorem}
\label{theorem:equilibrium}
For a functional $F:\mathbb{R}^d\rightarrow\mathbb{R}$, a index-$k$ saddle point $\boldsymbol{u}^*$ with inner product $\langle \cdot,\cdot\rangle_{L}$ is also a index-$k$ saddle point with inner product $\langle\cdot,\cdot\rangle_{T}$.
\end{theorem}
\begin{Proof}

By definition, $\boldsymbol{u}^*$ being an index-$k$ saddle point under $\langle \cdot,\cdot \rangle_L$ means
\[
\nabla_L F(\boldsymbol{u}^*) = 0 \quad \text{and} \quad \nabla^2_L F(\boldsymbol{u}^*) \text{ has exactly } k \text{ negative eigenvalues}.
\]
Since $\nabla_L F(\boldsymbol{u}^*) = \operatorname{L} \nabla F(\boldsymbol{u}^*)$, and $\ker(\operatorname{L}) = span\{\boldsymbol{1}\}$, we have
\[
\nabla_L F(\boldsymbol{u}^*) = 0 \;\Leftrightarrow\; \nabla F(\boldsymbol{u}^*) \in span\{\boldsymbol{1}\},
\]
i.e., $\nabla F(\boldsymbol{u}^*)$ is a constant multiple of $\boldsymbol{1}$. 
Similarly, for inner product $\langle \cdot,\cdot\rangle_T$:
    \[
     \nabla_T F(\boldsymbol{u}^*)=0 \Leftrightarrow  \nabla F(\boldsymbol{u^*})\in span\{\boldsymbol{1}\}.
    \]
    Therefore, $\nabla_LF(\boldsymbol{u})=0 \Leftrightarrow \nabla_TF(\boldsymbol{u})=0$. \par
    Since $\operatorname{T_0}$ is self-adjoint and positive definite, we assume that $\operatorname{T_0}=(\operatorname{T}_0^{\frac{1}{2}})^2$, where $\operatorname{T}_0^{\frac{1}{2}}$ is self-adjoint and positive definite. From the process of proof Lemma \ref{FA_decomposed},
    \[
        \operatorname{T}\nabla^2F(\boldsymbol{u})[\operatorname{P}_h \boldsymbol{v}]=\lambda\boldsymbol{v}\Leftrightarrow \operatorname{T}_0^{\frac{1}{2}}\operatorname{P}_h\nabla^2F(\boldsymbol{u})[\operatorname{P}_h\operatorname{T}_0^{\frac{1}{2}} \boldsymbol{w}]=\lambda \boldsymbol{w},
    \]
    \begin{sloppypar}
        where $\boldsymbol{w}=(\operatorname{T}_0^{\frac{1}{2}})^{-1}\boldsymbol{u}$. Due to the contract between $\operatorname{T}_0^{\frac{1}{2}}\operatorname{P}_h\nabla^2F(\boldsymbol{u})\operatorname{P}_h\operatorname{T}_0^{\frac{1}{2}}$ and $\sqrt{\operatorname{L}_0}\operatorname{P}_h\nabla^2F(\boldsymbol{u})\operatorname{P}_h\sqrt{\operatorname{L}_0}$, the Sylvester's law of inertia indicates that the number of positive and negative eigenvalues is the same.
    \end{sloppypar}

\end{Proof}

\section{Numerical experiment}
In this section, we construct solution landscapes under various parameters and compare the results of the two dynamical systems based on discrete projected $L^2$ inner product and the $H^{-1}$ inner product for the Ginzburg-Landau model with driving force. 
First, we introduce the rescaled Ginzburg–Landau free energy functional with driving force, and, via spatial discretization, obtain two dynamical systems: one with the projected discrete $L^2$ inner product with projection and the other with the $H^{-1}$ inner product.
Second, we numerically verify the theoretical properties established in Section 2, including convergence rates. 
Finally, we construct and compare the solution landscapes of the two dynamics under Neumann and periodic boundary conditions, demonstrating how the choice of inner product enriches the solution landscape.

\subsection{Saddle dynamics and the numerical discretization for the 
driven Ginzburg–Landau free energy}

The rescaled Ginzburg-Landau functional with driving force \cite{2007PhysRevLett} is shown as follows:

\begin{equation}
	\mathcal{F}(\phi)= \int_{\Omega} \left( \frac{\epsilon^2}{2} |\nabla \phi|^2 + \frac{1}{4} L^2(\phi^2 - 1)^2  -\frac{1}{4}\lambda L^2(3\phi-\phi^3)\right) dx,
	\label{DF_functional_t}
\end{equation}
where the domain is $\Omega=[0,1]\times[0,1]\subset \mathbb{R}^2$. Here, the parameter $\epsilon$ denotes the interfacial thickness parameter, $L$ is the rescaled domain size, and the parameter $\lambda$ represents the driving coefficient, which breaks the symmetry by forcing the evolution toward one phase.\par

We first discretize the domain $[0,1]\times[0,1]\subset \mathbb{R}^2$ of the functional \eqref{DF_functional_t} with a uniform Cartesian grid of $128\times 128$ points. 


Then, the saddle dynamics derived from the projected $L^2$ inner product and the $H^{-1}$ inner product are given respectively. The $L^2$ inner product saddle dynamics under projection (PSD) is
\begin{align}
    \begin{cases}
    	\dot{\boldsymbol{u}}=-\left(I-2\sum_{i=1}^{k}\boldsymbol{v}_{i}\boldsymbol{v}_i^{w,\top}\right)\operatorname{P}_h\nabla F(\boldsymbol{u}), \\
    	\dot{\boldsymbol{v}}_{i}=-\left(I-\boldsymbol{v}_{i}\boldsymbol{v}_i^{w,\top}-2\sum_{j=1}^{i-1}\boldsymbol{v}_{j}\boldsymbol{v}_j^{w,\top}\right)\operatorname{P}_h\nabla^2F(\boldsymbol{u})[\operatorname{P}_h\boldsymbol{v}_{i}],\ i=1,\ldots,k.
    \end{cases}
    \label{PSD}
\end{align}
The $H^{-1}$ inner product saddle dynamics ($H^{-1}$-SD) is
\begin{align}
    \begin{cases}                 \dot{\boldsymbol{u}}=\left(I+2\sum_{i=1}^{k}\boldsymbol{v}_{i}\boldsymbol{v}_i^{w,\top}\Delta_h^{-1}\right)\operatorname{\Delta_h}\nabla F(\boldsymbol{u}), \\
    \dot{\boldsymbol{v}}_{i}=\left(\Delta_h+\boldsymbol{v}_{i}\boldsymbol{v}_i^{w,\top}\operatorname{P}_h+2\sum_{j=1}^{i-1}\boldsymbol{v}_{j}\boldsymbol{v}_j^{w,\top}\operatorname{P}_h\right)\nabla^2F(\boldsymbol{u})[\operatorname{P}_h\boldsymbol{v}_{i}],\ i=1,\ldots,k,
    \end{cases}
    \label{H_PSD}
\end{align}
where $\boldsymbol{v}_i^{w,\top}=(\boldsymbol{w}\circ\boldsymbol{v}_i)^\top$. The gradient and the Hessian of $F$ are given by:
\[
\nabla F(\boldsymbol{u})=-\epsilon^2 \Delta_h \boldsymbol{u}+L^2(\boldsymbol{u}^3-\boldsymbol{u})-\frac{3}{4}L^2\lambda(1-\boldsymbol{u}^2),\]
and
\[
\nabla^2F(\boldsymbol{u})=-\epsilon^2\Delta_h+L^2(\frac{3}{2}\lambda \boldsymbol{u}-1+3\boldsymbol{u}^2).
\]

The weight vector $\boldsymbol{w}$ and the discrete Laplacian $\Delta_h$ depend on the chosen spatial grid and boundary conditions. The weight $\boldsymbol{w}$ for the homogeneous Neumann boundary condition is 
\begin{equation}
(\boldsymbol{w})_i=
\begin{cases}
\frac{1}{4},\ i=0,128,16256,16384,\\
\frac{1}{2},\ else,\\
1,i=130,131,\ ....,255,258,259,...383,...,16255
\end{cases},
\nonumber
\end{equation}
and $\boldsymbol{w}=(1,...,1)^\top\in \mathbb{R}^d$ for the periodic boundary condition.\par
We note that $\Delta_h: \mathbb{R}^d \to \mathbb{R}^d$ denotes the discrete Laplacian operator. For any grid function $\boldsymbol{u}\in\mathbb{R}^d$, $(\Delta_h \boldsymbol{u})_i$ under the homogeneous Neumann boundary condition is defined by the finite difference scheme of the five-point stencil, which is a square matrix
\begin{equation}
    \nonumber
    \Delta_h=\begin{bmatrix}
        B_1 & I & 0 & \cdots & 0 & 0 \\
        I & B_2 & I & \cdots & 0 & 0 \\
        0 & I & B_2 & \cdots & 0 & 0 \\
        \vdots & \vdots & \vdots & \ddots & \vdots & \vdots \\
        0 & 0 & 0 & \cdots & B_2 & I \\
        0 & 0 & 0 & \cdots & I & B_1
    \end{bmatrix}_{128 \times 128},
\end{equation}
where $I$ is an identical matrix and $B_1,B_2 \in \mathbb{R}^{128\times128}$,
\begin{equation}
B_1 =
\begin{bmatrix}
-2 & 1 & 0 & \cdots & 0 & 0 \\
1 & -3 & 1 & \cdots & 0 & 0 \\
0 & 1 & -3 & \cdots & 0 & 0 \\
\vdots & \vdots & \vdots & \ddots & \vdots & \vdots \\
0 & 0 & 0 & \cdots & -3 & 1 \\
0 & 0 & 0 & \cdots & 1 & -2
\end{bmatrix},\ B_2 =
\begin{bmatrix}
-3 & 1 & 0 & \cdots & 0 & 0 \\
1 & -4 & 1 & \cdots & 0 & 0 \\
0 & 1 & -4 & \cdots & 0 & 0 \\
\vdots & \vdots & \vdots & \ddots & \vdots & \vdots \\
0 & 0 & 0 & \cdots & -4 & 1 \\
0 & 0 & 0 & \cdots & 1 & -3
\end{bmatrix}.
\nonumber
\end{equation}
The operator $\Delta_h$ is semi-positive definite, self-adjoint, and its kernel $\ker(\Delta_h)=span\{\boldsymbol{1}\}$. From Lemma \ref{lem:A_property}, the operator $\Delta_h^{-1}:\mathbb{R}^d_0 \to \mathbb{R}^d_0$ is inverse. From Theorem \ref{theorem:linear} and Theorem \ref{theorem:equilibrium}, saddle points of these dynamics systems \eqref{PSD} and \eqref{H_PSD} can be identified. Under periodic boundary conditions, the operator $\Delta_h$ is defined analogously using the five-point stencil with periodic extension; it is also symmetric, semi-positive definite, and satisfies $\ker(\Delta_h)=span\{\boldsymbol{1}\}$. In both cases, $\Delta_h^{-1}:\mathbb{R}^d_0 \rightarrow \mathbb{R}^d_0$ exists, and saddle points of the dynamics can be identified via Theorems \ref{theorem:linear} and \ref{theorem:equilibrium}.

For saddle dynamics \eqref{PSD}, we employ the explicit Euler method with step size $\delta t$. The low stiffness of the PSD enables stable explicit time stepping with relatively large time steps, as will be verified numerically. \par 
For dynamics \eqref{H_PSD}, a semi-implicit scheme is applied to enhance numerical stability in the iteration of $\boldsymbol{u}$ and $\boldsymbol{v}_i$ $(1\leq i \leq k)$, particularly in treating the high-order derivative terms $-\epsilon^2\Delta_h^2 \boldsymbol{u}$ and $-\epsilon^2\Delta_h^2 \boldsymbol{v}_{i},1\leq i \leq k$:

\begin{equation}
    \begin{cases}
        \dfrac{\boldsymbol{u}_{n+1}-\boldsymbol{u}_{n}}{\delta t}=-\epsilon^2\Delta_h^2\boldsymbol{u}_{n+1} +N_1(\boldsymbol{v}_{1,n},...,\boldsymbol{v}_{k,n},\boldsymbol{u}_n), \\
        \dfrac{\boldsymbol{v}_{i,n+1}-\boldsymbol{v}_{i,n}}{\delta t} =-\epsilon^2 \Delta_h^2\boldsymbol{v}_{i,n+1} + N_2(\boldsymbol{v}_{1,n},...,\boldsymbol{v}_{i,n},\boldsymbol{u}_n),1\leq i \leq k.
    \end{cases}
    \label{DP-CH_d}
\end{equation}

Where $N_1(\boldsymbol{u}_n,\boldsymbol{v}_{1,n},...,\boldsymbol{v}_{k,n})=\epsilon^2V_1\Delta_h^2\boldsymbol{u}_n+L^2(I-V_1)(\boldsymbol{u}_n^3+\frac{3}{4}\lambda \boldsymbol{u}_n^2 -\boldsymbol{u}_n-\frac{3}{4}\lambda)$, $V_1=2\sum_{i=1}^k\boldsymbol{v}_i(\boldsymbol{w}\circ\boldsymbol{v}_i)^\top (-\Delta_h)^{-1}$, $N_2(\boldsymbol{u}_n,\boldsymbol{v}_{1,n},...,\boldsymbol{v}_{i,n})=\epsilon^2V_2^i\Delta_h^2 \boldsymbol{v}_{i,n}+L^2(I-V_2^i)(\frac{3}{2}\lambda \boldsymbol{u}_n-1+3\boldsymbol{u}_n^2)\boldsymbol{v}_{i,n}$, $V_2^i=\boldsymbol{v}_{i,n}(\boldsymbol{w}\circ\boldsymbol{v}_{i,n})^\top(-\Delta_h)^{-1}+2\sum_{j=1}^{i-1}\boldsymbol{v}_{j,n}(\boldsymbol{w}\circ\boldsymbol{v}_{j,n})^\top(-\Delta_h)^{-1}$. We note that the linear equations obtained by the semi-implicit method \eqref{DP-CH_d} are uniquely solvable because the operator 
$\operatorname{I}+\epsilon^2\Delta^2_h\delta t$ is positive definite.\newline
\subsection{Convergence tests}
To ensure accuracy in subsequent comparisons, we test the convergence of both numerical schemes separately. The parameters are fixed as $\epsilon^2 = 0.01$, $\lambda = 1.2$, and $L = 1$.

For system (22) with $k = 1$, we compute the solution from $T = 0$ to $T = 0.1$ using the initial condition
\[
\boldsymbol{u}^0 = -\frac{3}{4}\lambda + 0.01 \cdot \boldsymbol{v}_2^0,
\]
where $\boldsymbol{v}_1^0$ and $\boldsymbol{v}_2^0$ are the eigenvectors corresponding to the smallest eigenvalues of $\nabla^2F_h(-\frac{3}{4}\lambda)$. A benchmark solution is obtained using a very small time step $\delta t = 10^{-3}/2^{10}$, while numerical solutions are computed with $\delta t = 10^{-3}/2^i$ for $1 \leq i \leq 7$. Table \ref{PSD convergence} shows errors $e(\boldsymbol{u}) = \|\boldsymbol{u} - \boldsymbol{u}_{\text{ref}}\|_\infty$, $e(\boldsymbol{v}_1) = \|\boldsymbol{v}_1 - \boldsymbol{v}_{1,\text{ref}}\|_\infty$, and the corresponding convergence orders.

\begin{table}[H]
\centering
\caption{Time accuracy test in explicit Euler PSD format: The errors and convergence order between the numerical solution and the benchmark solution at different time steps, where $T =0.1$}
\label{PSD convergence}
\begin{tabular}{|c|cccccc|}
\hline
$10^3 \times\delta t$ & $\frac{1}{2}$ & $\frac{1}{4}$ & $\frac{1}{8}$ & $\frac{1}{16}$ & $\frac{1}{32}$ & $\frac{1}{64}$  \\
\hline
$e(\boldsymbol{u})$ & 4.629e-7 & 2.310e-7 & 1.151e-7 & 5.680e-8 & 2.795e-8 & 1.352e-8 \\
order & & 1.003 & 1.006 & 1.018 & 1.023 & 1.047 \\
$e(\boldsymbol{v}_1)$ & 1.6942e-5 & 8.4545e-6 & 4.210e-6 & 2.158e-6 & 1.062e-6 & 5.138e-7 \\
order & & 1.003 & 1.006 & 1.018 & 1.023 & 1.047 \\
\hline
\end{tabular}
\end{table}
A similar convergence test is performed for system \eqref{H_PSD} with the same parameters over time $T\in[0,0.01]$. We take the solution generated by $\delta t = 10^{-4}/2^{10}$ as the benchmark solution. Tables \ref{CH convergence} present the errors and the convergence rates for the schemes with the time step size being halved from $\delta t = 10^{-4}/2$ to $\delta t =10^{-4}/2^6$.
\begin{table}[H]
\centering
\caption{Time accuracy test in semi-implicit $H^{-1}$-SD format: The errors and convergence order between the numerical solution and the benchmark solution at different time steps, where $T =0.01$}
\label{CH convergence}
\begin{tabular}{|c|cccccc|}
\hline
$10^4 \times\delta t$ & $\frac{1}{2}$ & $\frac{1}{4}$ & $\frac{1}{8}$ & $\frac{1}{16}$ & $\frac{1}{32}$ & $\frac{1}{64}$ \\
\hline
$e(\boldsymbol{u})$ & 1.041e-6 & 5.203e-7 & 2.596e-7 & 1.293e-7 & 6.416e-8 & 3.157e-8  \\
order & & 1.001 & 1.002 & 1.006 & 1.011 & 1.023\\
$e(\boldsymbol{v}_1)$ & 2.253e-4 & 1.126e-4 & 5.621e-5 & 2.800e-5 & 1.399e-5 & 6.836e-6 \\
order & & 1.000 & 1.002 & 1.005 & 1.011 & 1.023 \\
\hline
\end{tabular}
\end{table}
As shown in Tables \ref{PSD convergence} and \ref{CH convergence}, both the PSD and $H^{-1}$-SD dynamics achieve first-order convergence in time, consistent with theoretical expectations for the explicit Euler and semi-implicit schemes, respectively.\\
\textbf{Remark 2}\\
In the accuracy test, we employ a discrete numerical scheme for the equation of  $\boldsymbol{v}_i$ of PSD \eqref{PSD} and $H^{-1}$-SD \eqref{H_PSD} rather than the locally optimal block preconditioned conjugate gradient (LOBPCG) method because the convergence accuracy of LOBPCG is time-independent and would disrupt the first-order convergence of errors in $\boldsymbol{u}$ and $\boldsymbol{v}_i$. Actually, we adopt LOBPCG to compute $\boldsymbol{v}_i$ $(1\leq i \leq k)$ in subsequent simulations for its efficiency, accuracy, and stability.

\subsection{Solution landscape of driven Ginzburg-Landau free energy}
For the driven Ginzburg-Landau functional \eqref{DF_functional_t}, there exist multiple critical points and their connections that may be different in different dynamics. Now we compare solution landscapes constructed by PSD \eqref{PSD} and $H^{-1}$-SD \eqref{H_PSD}. For this driven functional \eqref{DF_functional_t}, the function $\phi_0 = -\frac{3}{4}\lambda$ is a critical point. The eigenvalues of  $\nabla^2 F(\phi_0)$ under homogeneous Neumann boundary conditions are given by
$$
\lambda_{k,m} = \epsilon^2 \pi^2 (k^2 + m^2) + L^2 \left( \frac{9}{16} \lambda^2 - 1 \right), \ k, m \in \mathbb{N}^+.
$$
The eigenvalues of $\nabla^2 F(\phi_0)$ under periodic boundary conditions are given by
$$
\lambda_{k,m} = 4 \epsilon^2 \pi^2 (k^2 + m^2) + L^2 \left( \frac{9}{16} \lambda^2 - 1 \right), \quad k, m \in \mathbb{Z}
\setminus \{0\}.
$$
The index of the saddle point $\phi_0$ depends on the parameters $\lambda$, $\epsilon$ and domain size $L$.

We will apply downward research \cite{2020Construction} to systematically construct solution landscapes. The method starts from an index-$k$ saddle point $\boldsymbol{u}^*$ as the parent state to locate an index-$s$ saddle point with $s < k$. Let $\boldsymbol{v}_1^*, \ldots, \boldsymbol{v}_k^*$ be the eigenvectors of $\boldsymbol{u}^*$. We introduce a small perturbation $\boldsymbol{u}^* + \epsilon\ \boldsymbol{v}^*$, where $\boldsymbol{v}^*$ is a linear combination of $\boldsymbol{v}_1^*, \ldots, \boldsymbol{v}_k^*$ and $\epsilon$ is a small constant. Taking $\boldsymbol{v}_1^*, \ldots, \boldsymbol{v}_s^*$ as the initial directions for $\boldsymbol{v}$, the saddle points are identified using the dynamics \eqref{PSD} and \eqref{H_PSD}. Therefore, we need to find a critical point as the parent state.

\subsubsection{Neumann boundary condition}
\begin{figure}[ht]
	\centering
    \includegraphics[scale=0.6]{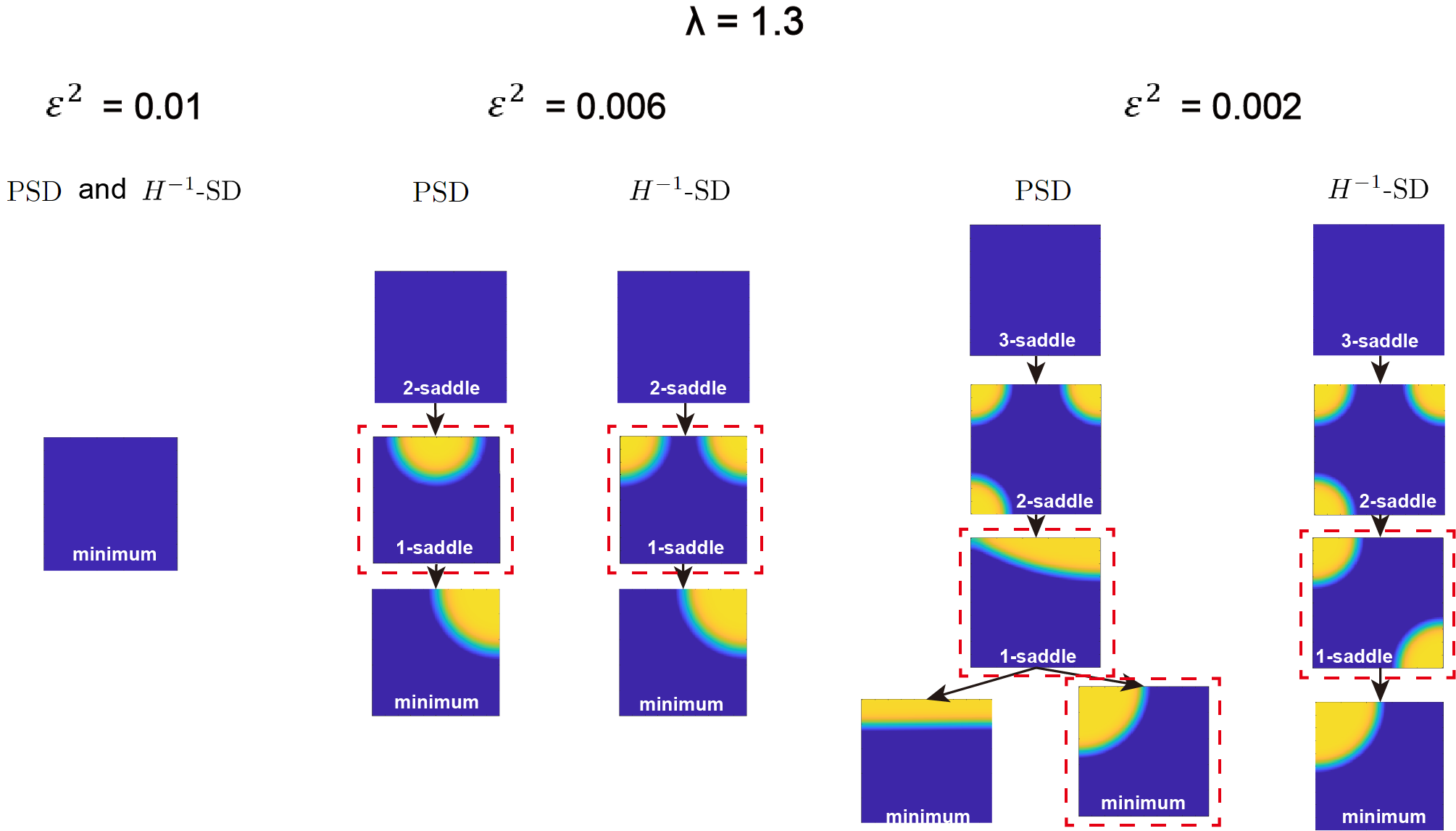}
    \caption{The solution landscape derived from PSD and $H^{-1}$-SD with $\lambda=1.3$ , $\epsilon^2=0.01,0.006,0.002$. The different solutions are framed with red square.}
    \label{lambda1_3N}
\end{figure}

\begin{figure}[h]
	\centering
    \includegraphics[scale=0.55]{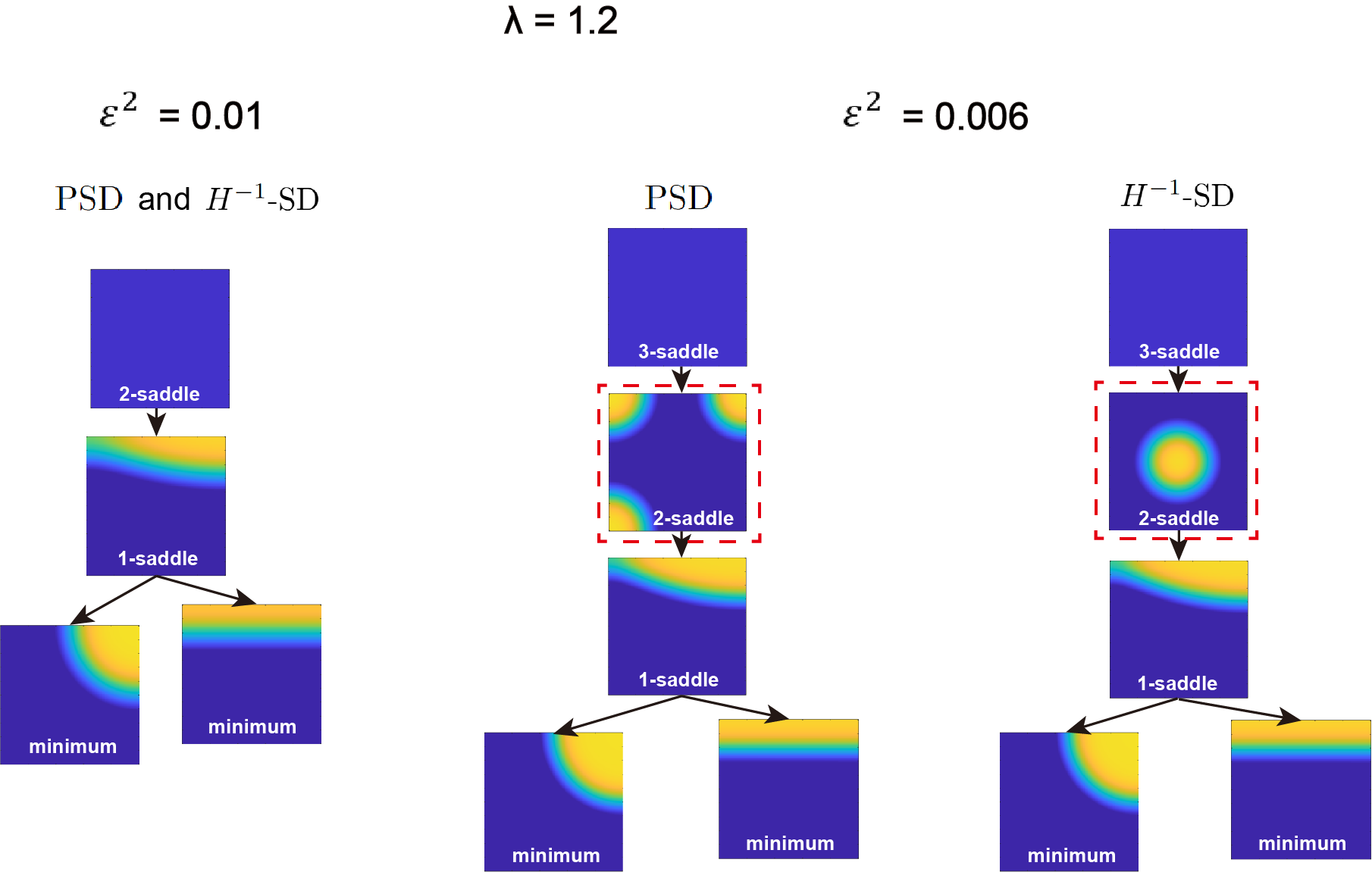}
    \caption{The solution landscape derived from PSD and $H^{-1}$-SD with $\lambda=1.2$, $\epsilon^2=0.01,0.006$. The different solutions are framed with red square.}
    \label{lambda1_2N}
\end{figure}

We apply \eqref{PSD} and \eqref{H_PSD} respectively to construct the solution landscapes with different parameters under the Neumann boundary condition. We first select $\lambda=1.3$ with $\epsilon^2=0.01,0.006,0.002$ and $\lambda=1.2$ with $\epsilon^2=0.01,0.006$. The solution landscapes with values $\phi_0=-\frac{3}{4}\lambda$ as parent states are then constructed in Figure \ref{lambda1_3N} and Figure \ref{lambda1_2N}. These parameter choices are intended to vary the parent state, allowing exploration of distinct solution landscapes with different topological structures by downward searching from the parent state.\par

In general, a larger $\lambda$ or smaller $\epsilon$ increases the index of the parent state $\phi_0$, leading to a more complex solution landscape. Moreover, the two saddle dynamics methods are likely to converge to different saddle points.

From the Figure \ref{lambda1_3N}, it is observed that the homogeneous solution is the minimum point with $\lambda=1.3$ and $\epsilon^2=0.01$. 
When $\epsilon^2$ decreases to $0.006$, the minimum becomes unstable and transforms to an index-$2$ saddle point of a quarter‑circle structure. At the index-$1$ saddle level, the two methods yield different solutions: PSD gives a semi‑circular structure, while $H^{-1}$-SD produces a two‑quarter‑circle structure. These two saddles are not equivalent under the symmetry of the domain, reflecting the influence of the inner product on the search path.
When $\epsilon^2$ decreases to $0.002$, the solution landscape becomes more diverse. The PSD method identifies two distinct minimizers: the layer solution and the quarter‑circle structure, whereas the layer solution is absent in the $H^{-1}$-SD landscape. At the index-$1$ saddle level, PSD yields a knife‑shaped structure, while $H^{-1}$-SD produces a two‑quarter‑circle configuration.

From Figure \ref{lambda1_2N}, we consider the parameter $\lambda = 1.2$ with interface parameters $\epsilon^2 = 0.01$ and $0.006$. For $\epsilon^2 = 0.01$, the two solution landscapes from PSD and $H^{-1}$-SD are same. When $\epsilon^2$ is reduced to $0.006$, the solution landscapes differ at the index-$2$ saddle point: PSD saddle dynamics yields three quarter-circle structures, while the $H^{-1}$ inner produces a circle at the center.

\begin{figure}[ht]
	\centering
    \includegraphics[scale=0.6]{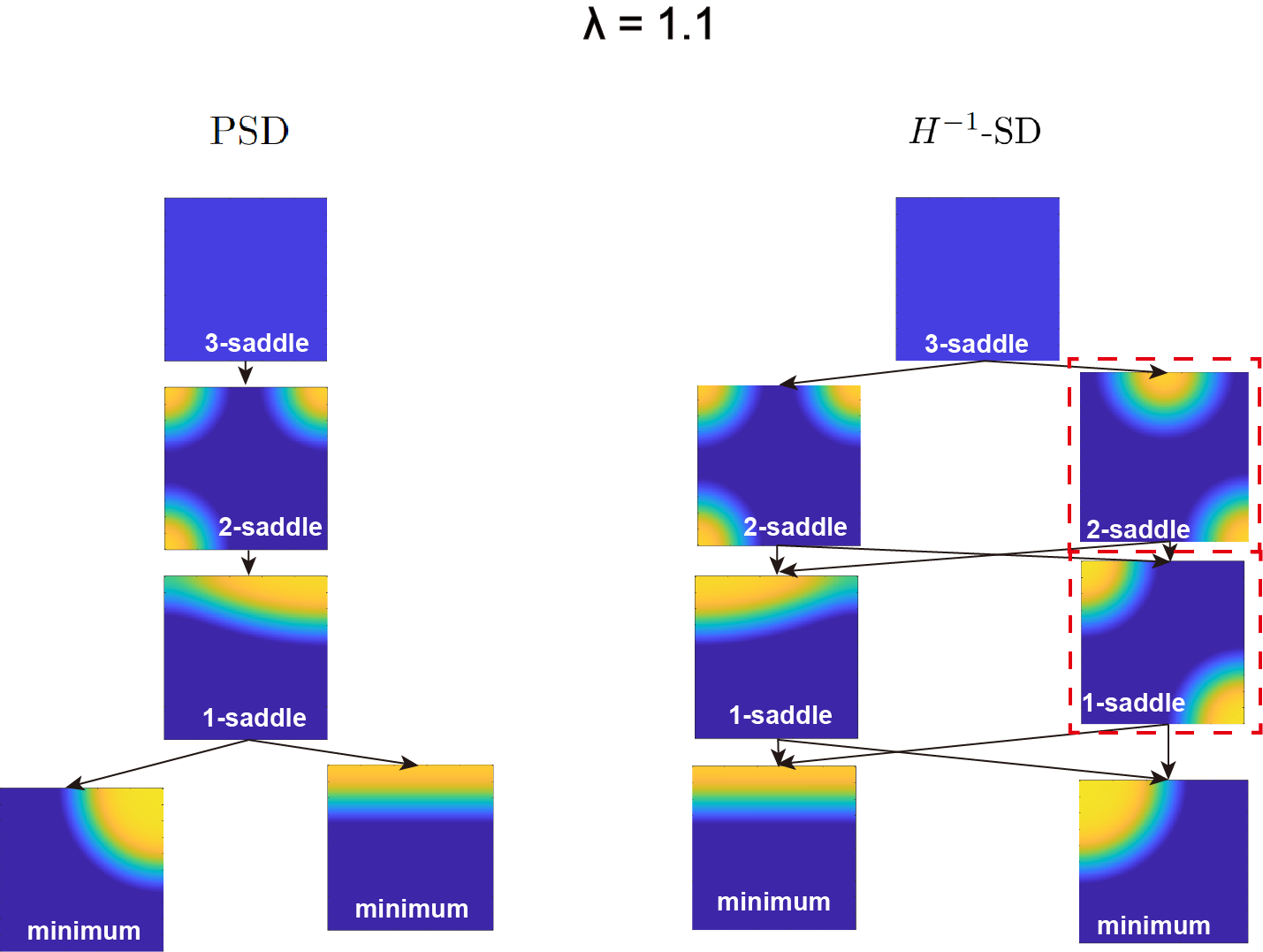}
    \caption{The solution landscape derived from PSD and $H^{-1}$-SD with $\lambda=1.1$,$\epsilon^2=0.01$. The different solutions are framed with red square.}
    \label{lambda1.1N}
\end{figure}

We then set $\lambda = 1.1$, for which the homogeneous solution $\phi_0$ becomes index-$3$. As shown in Figure~\ref{lambda1.1N}, the $H^{-1}$ inner product saddle dynamics yields additional solutions beyond those obtainable with PSD, including two index-$2$ solutions: a "half-circle-and-quarter" solution and a "two-quarter" solution.

\begin{figure}[h]
	\centering
    \includegraphics[scale=0.6]{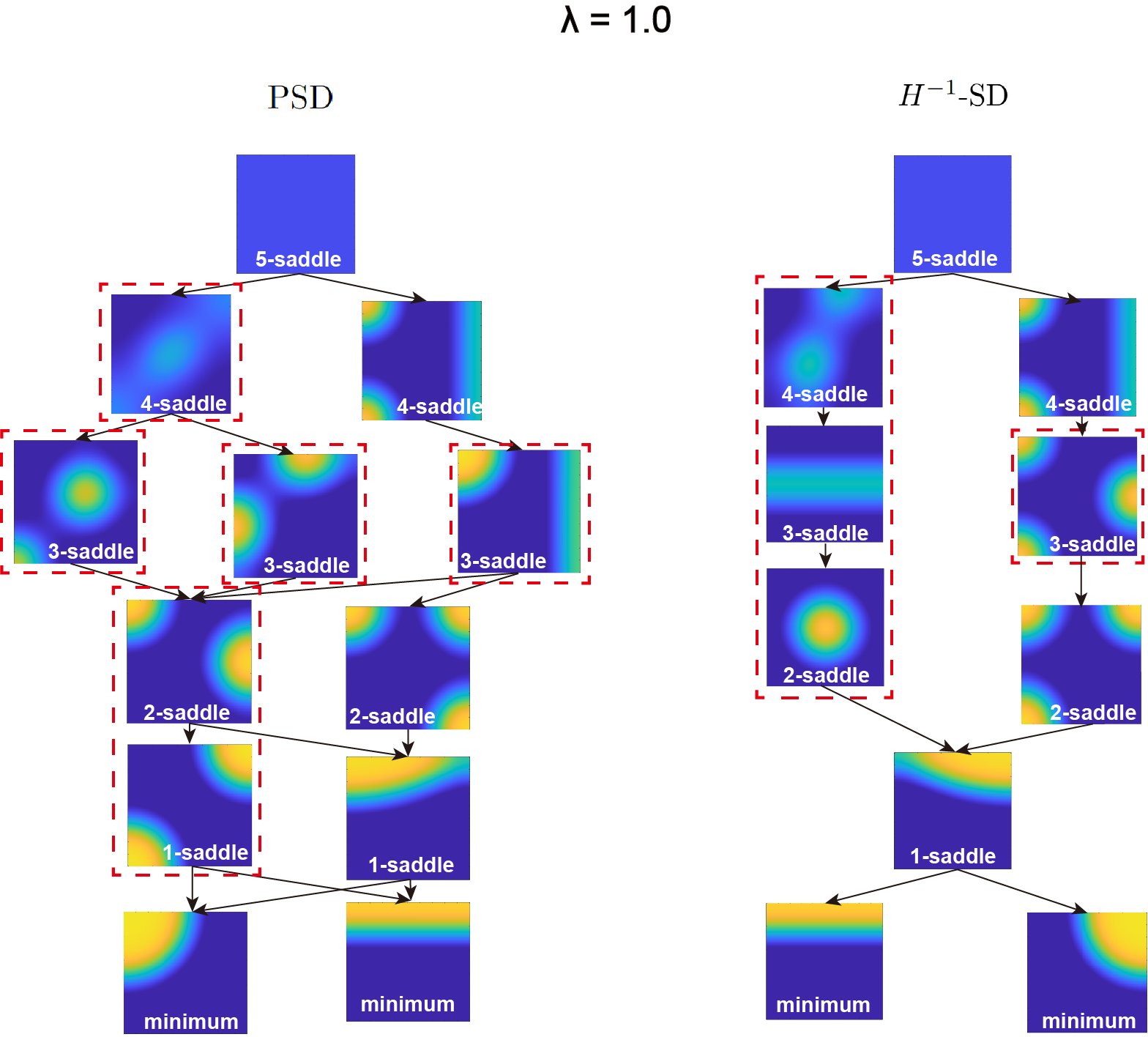}
    \caption{The solution landscape derived from PSD and $H^{-1}$-SD with $\lambda=1.0,\epsilon=0.01$. The different solutions are framed with red square.}
    \label{lambda1.0N}
\end{figure}

Further reduction of $\lambda$ to $1.0$ produces a richer solution landscape, as shown in Figure \ref{lambda1.0N}. The solution landscapes obtained with the inner products of PSD and $H^{-1}$ diverge further, and the saddle points of index-$1$ to index-$4$ are not identical.\par
\newpage
In conclusion, although the parent state of the solution landscape is the same, the saddle points and their connections differ between the inner product formulations of PSD and $H^{-1}$, due to their distinct evolution paths. Although PSD is more efficient, $H^{-1}$-SD can identify additional saddle points and relationships that PSD does not capture easily. In particular, the same saddle points of both methods are consistent with Theorem \ref{theorem:equilibrium}.
\subsubsection{Periodical boundary condition}
\begin{figure}[ht]
	\centering
    \includegraphics[scale=0.5]{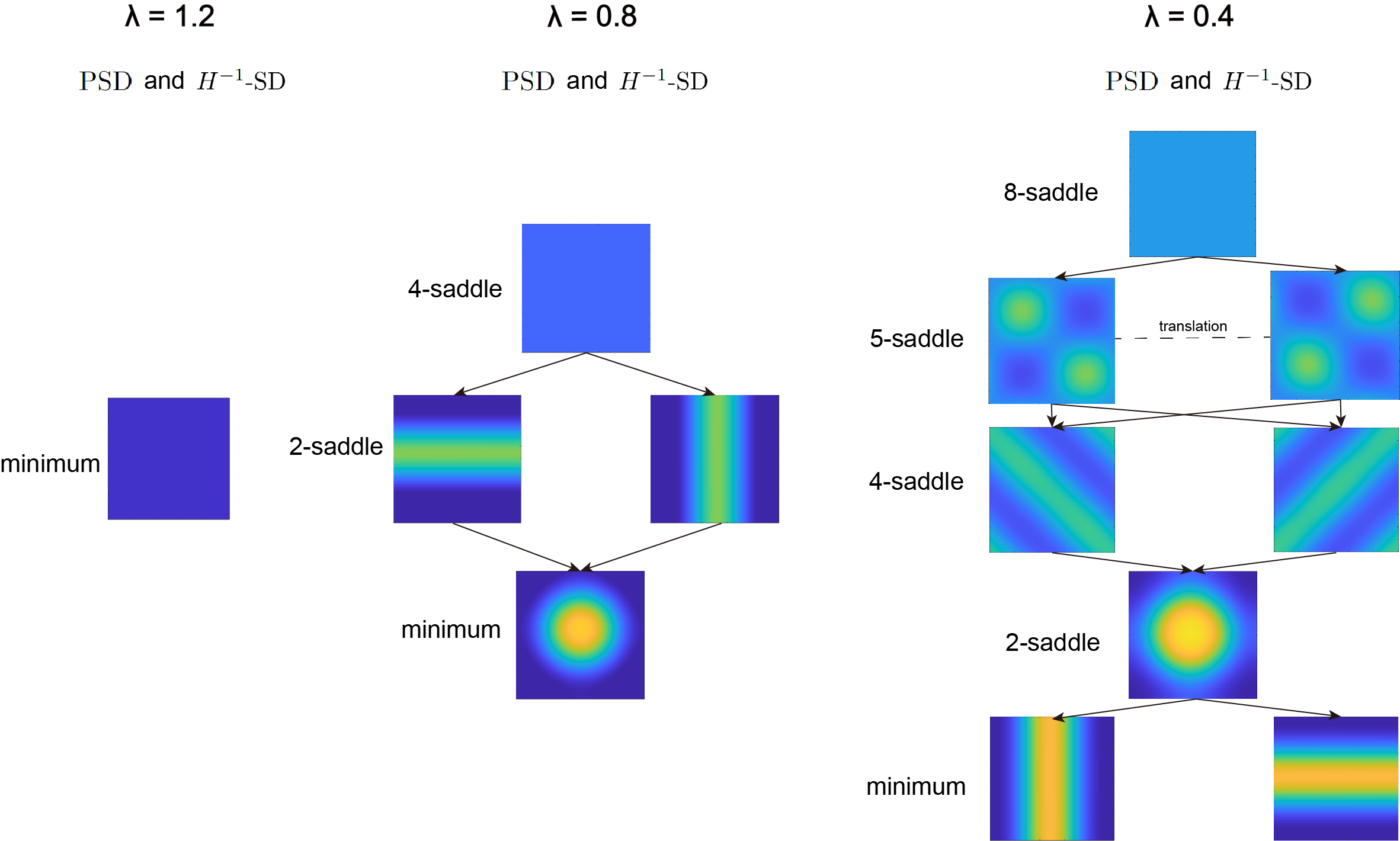}
    \caption{The solution landscape derived from PSD and $H^{-1}$-SD with $\lambda=1.2,0.8,0.4$.}
    \label{lambdaP_lambda}
\end{figure}

\begin{figure}[h]
	\centering
    \includegraphics[scale=0.8]{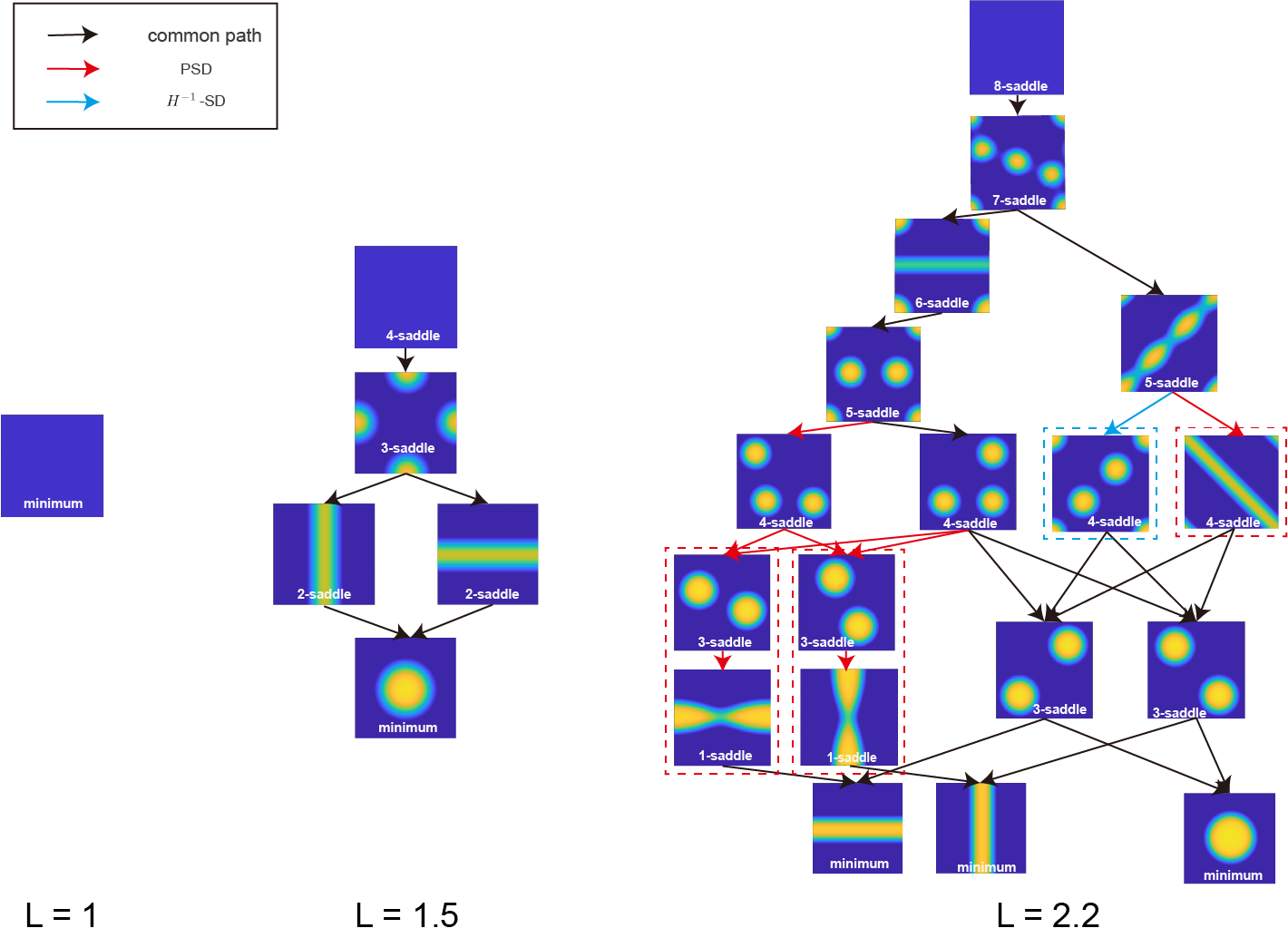}
    \caption{The solution landscape derived from PSD and $H^{-1}$-SD with $L=1,1.5,2.2$. Black arrows denote connections that both the $H^{-1}$ and PSD methods can identify. Red arrows represent connections found only by the PSD method, typically requiring a small initial perturbation. Cyan arrows indicate connections found only by the $H^{-1}$-SD method.}
    \label{lambdaP_L}
\end{figure}

To examine the impact of boundary conditions on the solution landscapes, we apply the same method to construct the solution of functional \eqref{DF_functional_t} with periodic boundary condition on $[0,1]\times[0,1]$. Firstly, we take the interface parameter $\epsilon^2=0.01$ and the driven parameter $\lambda =1.2,0.8,0.4$. Starting from the homogeneous solution $\phi^0=-\frac{3}{4}\lambda$, the two distinct dynamical systems derive the same solution landscape, as shown in Figure \ref{lambdaP_lambda}.
Secondly, we employ both the $H^{-1}$-SD and the PSD method to construct the solution landscape. To investigate the effect of the domain size $L$, we perform numerical experiments with $L = 1, 1.5, 2.2$, keeping the parameters fixed at $\epsilon^2 = 0.01$ and $\lambda = 1.2$. The resulting solution landscapes are presented in Figure \ref{lambdaP_L}. In the figure, connections between solutions are indicated by arrows of different colors. Black arrows denote connections that both the $H^{-1}$ and PSD methods can identify. Red arrows represent connections found only by the PSD method, typically requiring a small initial perturbation. Cyan arrows indicate connections found only by the $H^{-1}$-SD method.\par
For the domain size $L=1$ see first column of Figure \ref{lambdaP_L}, the homogeneous solution $\phi^0=-\frac{3}{4}$ is a minimum. When $L$ increases to $1.5$, $\phi^0 = -\frac{3}{4}$ becomes an index-$4$ saddle point. At this scale, the solution landscapes generated by the two methods are identical, yielding the same set of solutions and connecting pathways. The landscape becomes significantly more complex at $L = 2.2$, where $\phi^0 = -\frac{3}{4}$ is now an index-$8$ saddle point. A key observation is that two distinct index-$4$ saddle points—one characterized by a linear three-circle pattern and the other by a strip pattern—can be reached from a common index-$5$ strip-shaped saddle point. Notably, the linear-pattern index-$4$ saddle is found using $H^{-1}$-SD, whereas the strip-pattern index-$4$ saddle is located via the PSD method. Similarly, these index-$3$ and index-$1$ saddle points with red frame can be located with PSD method.

\section{Conclusion}
This work presents the Conservative Generalized Inner Product Saddle Dynamics (CGiSD) as a unified mathematical framework for saddle-point search under mass-conservation constraints. By embedding the conservation constraint into a generalized inner product induced by a self-adjoint, positive semi-definite operator whose kernel is the constant subspace, the framework naturally encompasses both projection saddle dynamics (PSD) and diffusion-based dynamics (e.g., the $H^{-1}$ inner product dynamics) as special cases. We establish the well-posedness of CGiSD, prove its linear stability, and demonstrate the equivalence of saddle points under different inner products, thereby ensuring consistency in steady-state predictions.
Numerically, we implement CGiSD with two discrete inner products $L^2$ with projection and $H^{-1}$ applied to the Ginzburg–Landau functional with a driving force, serving as a prototype model for symmetry-breaking phase transitions. Through extensive simulations under both Neumann and periodic boundary conditions, we systematically compare the solution landscapes obtained by the two dynamics and apply them to enrich the solution landscapes.\par
The results of the article raise several intriguing questions. First, the convergence analysis of discrete CGiSD is important open problem; the conservative generalized inner product can be viewed as a preconditioner in saddle dynamics, potentially leading to faster convergence. Second, high-order time discretization methods can be applied to CGiSD, such as Runge–Kutta, Adams–Bashforth, and Adams–Moulton. A key direction for future work is to rigorously establish convergence frameworks and extend CGiSD to a broader class of models and inner products.




\section*{Acknowledgments}
G. Ji is partially supported by the National Natural Science Foundation of China (Grant No.12471363). Z. Xu is partially supported by the National Natural Science Foundation of China (Grant No.12201021). The authors thank Dr. J. Yin for valuable comments and suggestions.







\bibliographystyle{cas-model2-names}
\bibliography{cas-refs}



\end{document}